\RequirePackage{fix-cm}
\documentclass[smallextended]{svjour3}       
\smartqed  
\usepackage{graphicx}
\usepackage{mathptmx}      
\usepackage{amsmath,amssymb,amsfonts}
\usepackage{color}
\usepackage{hyperref}
\usepackage{diagbox}
\usepackage{subcaption}
\usepackage[ruled,vlined,linesnumbered]{algorithm2e}
\usepackage{float}
\allowdisplaybreaks

\usepackage{etoolbox}
\newcommand*{\affaddr}[1]{#1} 
\newcommand*{\affmark}[1][*]{\textsuperscript{#1}}

\newcommand{\cref}[1]{(\ref{#1})}

\journalname{Journal of Scientific Computing}
\begin{document}

\title{A deep shotgun method for solving high-dimensional parabolic partial differential equations}
\titlerunning{Deep shotgun method}        

\author{Wenjun Xu\affmark[1,2]         \and
        Wenzhong Zhang\affmark[2,3] 
}


\institute{Wenjun Xu \at
              \email{wjxu@mail.ustc.edu.cn}           
           \and
           Wenzhong Zhang \at
              \email{wenzhong@ustc.edu.cn}
              \\
              \\
\affaddr{\affmark[1] School of Mathematical Sciences, University of Science and Technology of China, Hefei 230026, Anhui, P. R. China} \\
\affaddr{\affmark[2] Suzhou Institute for Advanced Research, University of Science and Technology of China, Suzhou 215123, Jiangsu, P. R. China} \\
\affaddr{\affmark[3] Key Laboratory of the Ministry of Education for Mathematical Foundations and Applications of Digital Technology, University of Science and Technology of China, Hefei 230027, Anhui, P. R. China}
}

\date{Received: date / Accepted: date}

\maketitle

\begin{abstract}
Recent advances in deep learning makes solving parabolic partial differential equations (PDEs) in high dimensional spaces possible via forward-backward stochastic differential equation (FBSDE) formulations.
The implementation of most existing methods requires simulating multiple trajectories of stochastic processes with a small step size of time discretization to ensure accuracy, hence having limited performance, especially when solving on a large time interval.
To address such issue, we propose a deep ``shotgun method'' that does not exploit full trajectories, but only utilizes the data distribution of them.
Numerical results including examples with dimensionality up to $10000$ demonstrate the competitiveness of the proposed shotgun method in both performance and accuracy.
\end{abstract}

\section{Introduction}\label{sec_intro}
Solving partial differential equations (PDEs) numerically in high-dimensional spaces has been known as a challenging task because of the curse of dimensionality (CoD), namely, the number of unknowns in the numerical solver grows exponentially as the dimensionality increases.
For high-dimensional parabolic PDEs such as the Black--Scholes equation for option pricing and the Hamilton--Jacobi--Bellman equation originated in stochastic optimal control problems, recent approaches such as DeepBSDE \cite{r1}, FBSNNs \cite{r2} and derived work \cite{r9}, the diffusion loss schemes \cite{r10,DiffusionLoss24}, DeepMartNet \cite{r11,cai2024soc}, etc. take advantage of deep learning, where the deep neural networks (DNNs) are capable of approximating high-dimensional functions using a feasible number of trainable variables.
Note that there's no intention for a thorough enumeration.
Unlike in the physics-informed neural network (PINN) \cite{r3} method that minimizes the mean squared residual of a PDE at spatially sampled positions, which is arguably the most prevailing general-purpose deep learning method for solving PDEs numerically, the loss function and the data collection of the above-mentioned methods are derived from the Pardoux--Peng theory \cite{r8}  of nonlinear Feynman--Kac formula.
The Pardoux--Peng theory bridges a parabolic PDE with a forward-backward stochastic differential equation (FBSDE) system.
By matching the PDE formulation and the FBSDE formulation in different manners, different loss function schemes have been developed in these ``deep FBSDE'' numerical methods.
Meanwhile, the data collection is similar across all these deep FBSDE methods: along discretized trajectories of the stochastic processes in the FBSDE formulation, towards a prescribed stopping time.

There are two common benefits of the deep FBSDE methods compared to PINN for numerically solving high-dimensional parabolic PDEs.
First, they effectively prevent explicit evaluation of the second-order partial derivatives in the PDE formulation, which, within typical deep learning frameworks, requires the full $d \times d$ Hessian matrix and takes $\mathcal{O}(d^2)$ computational cost when the spatial dimensionality is $d \gg 1$.
Second, the stochastic process trajectories in the FBSDE formulation are usually meaningful in practical problems, hence providing a helpful suggestion to the data distribution used for training.
However, when simulating the trajectories in the deep FBSDE methods in high-dimensional space, the Euler--Maruyama method is the only discretization approach of decent efficiency, but its accuracy is limited to half order with respect to time discretization step size.
Thus, existing deep FBSDE methods have to apply large discretization step numbers for better resulting accuracy.
In some of the implementations, the outputs of DNNs along the trajectories are accumulated before backward propagation, which further slows down the training process.
As a result, the capability of solving in a long time range is limited.
To address such issue, we propose a deep ``shotgun method'' that utilizes the data distribution of trajectories instead of the full-length version considering all the time steps.
This is made possible by the observation that each single-step discretization of the FBSDE system on a time interval with length $0 < \Delta t \ll 1$ is indeed an estimator of the PDE residual.
The estimate gets improved as $\Delta t \to 0^+$, and with multiple independently sampled antithetic variates.
The sample points are still picked from stochastic process trajectories, but using a much larger step size for simulation is sufficient in our numerical tests.
The resulting method behaves like shot grouping on the selected sample points from trajectories, hence given the same name as in the shotgun sequencing of DNA strands\cite{r5,r6}.

The rest of this paper is organized as follows.
In Section \ref{sec_related_work}, we take a brief review of a few deep FBSDE methods with our concerns.
Then, we introduce the deep shotgun method in Section \ref{sec_method}.
The numerical tests are presented in Section \ref{sec_result}, with succeeding discussion in Section \ref{sec_discussion}.
Some final concluding remarks are given in Section \ref{sec_conclusion}.

\section{A review of deep FBSDE methods}\label{sec_related_work}
Before introducing the shotgun method, we will first briefly review selected existing deep FBSDE methods that utilizes a FBSDE formulation to solve parabolic PDEs with deep learning.
\subsection{The parabolic PDE and the FBSDE formulation}\label{sec_pardoux_peng}
We study $d$-dimensional parabolic PDEs of $u(t, x)$ in the form
\begin{align}\label{pde}
\begin{split}
\mathcal{L}[u]: = \partial_t u + \frac{1}{2}\mathrm{Tr}[\sigma \sigma^T \nabla \nabla u] + \mu \cdot \nabla u &= \phi,\quad t\in [0, T], \quad x \in \mathbb{R}^d, \\
u(T,x) &= g(x),\quad x \in \mathbb{R}^d
\end{split}
\end{align}
where $\mu$, $\sigma$, $\phi$ are functions of $t, x, u(t, x), \nabla u(t, x)$, $\cdots$ that map to $\mathbb{R}^d$, $\mathbb{R}^{d \times d}$ and $\mathbb{R}$, respectively, and $g:\mathbb{R}^d \to \mathbb{R}$ presents the terminal condition values at time $T$.
For the rest of this paper, we use $\nabla$ for the gradient with respect to spatial coordinates $x$.
The goal of solving the PDE is usually finding an initial value at given $x_0 \in \mathbb{R}^d$, i.e. finding $u(0, x_0)$.

Under the Pardoux--Peng theory \cite{r4,r7}, the PDE \cref{pde} is associated with the following FBSDE system
\begin{align}\label{fbsde}
    \begin{split}
        d X_t&= \mu dt +\sigma d W _t, \quad t\in[0,T], \\
        X_0 &= x_0, \\
        d Y_t&=\phi dt+{Z_t}^T\sigma d W_t, \quad t\in[0,T], \\
        Y_T &= g(X_T),
    \end{split}
\end{align}
where $\{W_t: 0 \le t \le T\}$ is a $d$-dimensional standard Brownian motion, and the variables of functions $\mu$, $\sigma$ and $\phi$ consist of $t, X_t, Y_t, Z_t, \cdots$ instead of $t, x, u(t, x), \nabla u(t, x)$, $ \cdots$, respectively.
Indeed, the relationship between the stochastic process $X_t$, $Y_t$ and $Z_t$ in the FBSDE is given by the solution to the PDE as
\begin{equation}\label{link_PDE_FBSDE}
    Y_t = u(t, X_t), \quad Z_t = \nabla u(t, X_t).
\end{equation}
This can be verified using the Itô's formula
\begin{equation*}
du(t, X_t) = \left(\partial_t u + \frac{1}{2} \sigma^T \nabla \nabla u \sigma \right) dt + \nabla u^T (\mu dt + \sigma dW_t) = \mathcal{L}[u] dt + \nabla u^T \sigma dW_t.
\end{equation*}
With certain regularity conditions imposed on the functions $\mu, \sigma$ and $\phi$, the FBSDE system has unique adapted solution.

Numerically, the FBSDE system is usually discretized using the Euler--Maruyama method.
Here, we consider the uniform time discretization
\begin{equation*}
    0 = t_0 < t_1 < \cdots < t_N = T.
\end{equation*}
The SDE of $X_t$ and $Y_t$ are discretized as
\begin{align}\label{eq_FBSDE_discretization}
\begin{split}
    X_{n+1} &= X_n + \mu(t_n, X_n, Y_n, Z_n, \cdots) \Delta t_n^1 + \sigma(t_n, X_n, Y_n, Z_n, \cdots) \Delta W_n, \\
    Y_{n+1} &= Y_n + \phi(t_n, X_n, Y_n, Z_n, \cdots) \Delta t_n^1 + Z_n^T \sigma(t_n, X_n, Y_n, Z_n, \cdots) \Delta W_n,
\end{split}
\end{align}
respectively, where $\Delta t_n^1 = t_{n+1} - t_{n}$ is the increment of time, and $\Delta W_n = W_{n+1} - W_n$ is the increment of the Brownian motion on $[t_n, t_{n+1}]$.
The derived $X_n$, $Y_n$ and $Z_n$ are approximations of $X_{t_n}$, $Y_{t_n}$ and $Z_{t_n}$ from the continuous-time stochastic processes, respectively.

\subsection{Deep FBSDE methods}\label{sec_review_deep_fbsde}
Here we enumerate a few deep FBSDE methods and discuss how the FBSDE formulation facilitates deep learning approaches.
\subsubsection{Deep BSDE}
The Deep BSDE method \cite{r1} trains a DNN that accepts an instance of the Brownian motion $W_t$ as the input, and takes $Y_N$ as the output.
The initial values at $X_0 = x_0$ are demanded, which are trainable variables $Y_0 = u(0, x_0)$ and $Z_0 = \nabla u(0, x_0)$ in the network.
The value of $Y_N$, which is an approximation to $u(t_N, X_N)$, is derived by simulating the $X_t$ and $Y_t$ trajectories, where in the discretization, sub-networks $f^{(n+1)}_\theta \approx \nabla u(t_{n+1}, \cdot)$ are applied to take the place of the discretized unknown process $Z_{n+1}$.
The loss function is constructed to minimize a Monte Carlo approximation of
\begin{equation}
    \mathbb{E}[\lVert Y_N - g(X_N) \rVert^2].
\end{equation}
The Deep BSDE method has been shown to give convergent numerical results for various high-dimensional parabolic equations \cite{r1}, and a posteriori estimate suggests strong convergence of half order \cite{r12} with respect to time discretization.
However, a drawback of our concern is that the number of parameters in the entire network increases linearly with respect to the number $N$ of time discretization steps.

\subsubsection{FBSNNs and derived work}\label{FBSNNs}
The FBSNNs method \cite{r2} trains a DNN $u_{\theta} (t, x)$ that directly approximates the solution to the PDE \cref{pde} on the entire time interval $[0, T]$.
The method makes use of the two distinct representations of $Y_n = u(t_n, X_n)$ from the PDE and from the FBSDE formulations.
Namely, in the discretization \cref{eq_FBSDE_discretization} of $Y_t$, FBSNNs uses the DNN outputs $u_\theta(t_n, X_n)$ and $u_\theta(t_{n+1}, X_{n+1})$ to replace $Y_n$ and $Y_{n+1}$, respectively, to form an identity that is supposed to hold as $\Delta t_n \to 0$, provided $u_\theta$ is precisely the demanded solution.
A squared residual loss is formed along the $Y_t$ trajectory, together with the matching of terminal values $u_\theta(t_N, X_N) = g(X_N)$ and $\nabla u_\theta(t_N, X_N) = \nabla g(X_N)$.

In this paper, we formulate the FBSNNs method for one instance of $W_t$ trajectory as follows.
First, for the demanded initial value at $x_0$, let $X_0 = x_0$, $Y_0 = \hat{Y}_0 = u_\theta(0, X_0)$ and $Z_0 = \nabla u_\theta(0, X_0)$ with a random initialization of the DNN $u_\theta$.
Second, the trajectories are simulated with an iterative scheme
\begin{align}\label{eq_fbsnn_scheme}
\begin{split}
    X_{n+1} &= \mu(t_n, X_n, Y_n, Z_n, \cdots) \Delta t_n^1 + \sigma(t_n, X_n, Y_n, Z_n, \cdots) \Delta W_n, \\
    Y_{n+1} &= u_\theta(t_{n+1}, X_{n+1}), \\
    \hat{Y}_{n+1} &= Y_{n} + \phi(t_n, X_n, Y_n, Z_n, \cdots) \Delta t_n^1 + Z_n^T \sigma(t_n, X_n, Y_n, Z_n, \cdots) \Delta W_n, \\
    Z_{n+1} &= \nabla u_\theta(t_{n+1}, X_{n+1}),
\end{split}
\end{align}
until $t_N = T$.
Finally, the loss function
\begin{equation}
    \sum_{n=1}^{N} \lVert Y_n - \hat{Y}_n \rVert^2 + \lVert Y_N - g(X_N) \rVert^2 + \lVert Z_N - \nabla g(X_N) \rVert^2
\end{equation}
is minimized to train the DNN $u_\theta$.
The loss function used by FBSNNs consists of both terms from the terminal conditions at $t=T$, and terms accumulated from the trajectory.
The difference between $Y_{n+1} - \hat{Y}_{n+1}$ can be interpreted as a discretization of the left-hand side of the SDE
\begin{equation}
dY_t - \phi(t, X_t, Y_t, Z_t, \cdots) dt - Z_t^T \sigma(t, X_t, Y_t, Z_t, \cdots) dW_t = 0,
\end{equation}
which holds if $u_\theta$ is the solution to the PDE.

A derived work \cite{r9} claims that the above stochastic process $\hat{Y}_n$ may not be a proper discretization of a SDE, so an alternative approach is proposed by using
\begin{equation}
    \hat{Y}_{n+1} = \hat{Y}_{n} + \phi(t_n, X_n, Y_n, Z_n, \cdots) \Delta t_n^1 + Z_n^T \sigma(t_n, X_n, Y_n, Z_n, \cdots) \Delta W_n
\end{equation}
instead of the one used in \cref{eq_fbsnn_scheme}.
Thus, $Y_n$ is a discretization of the SDE $dY_t = u_\theta(t, X_t) dt$, and $\hat{Y}_n$ becomes a discretization of the SDE
\begin{equation}
d\hat{Y}_t = \phi(t, X_t, Y_t, Z_t, \cdots) dt + Z_t^T \sigma(t, X_t, Y_t, Z_t, \cdots) dW_t.    
\end{equation}
We refer to this approach (Scheme 2 from \cite{r9}) as a ``SDE matching'' method in this paper, as it attempts to match $dY_t$ and $d \hat{Y}_t$ along trajectories.
SDE matching achieves half order of accuracy with respect to time discretization in numerical tests.

Both FBSNNs \cite{r2} and its derived work \cite{r9} use a single DNN regardless of the time discretization.
However, we emphasize that they simulate complete trajectories.
Moreover, the evaluation of the term $\lVert Y_n - \hat{Y}_n \rVert^2$ in the loss function utilizes the output values of the DNN $u_\theta$ at all previous time steps, so the backward propagation of the loss function has at least $\mathcal{O}(N^2)$ computational cost.
When a uniform time discretization is applied, i.e. when each $\Delta t_n^1 = T / N$, a large $N$ is also demanded to reduce the size of $\Delta t_n^1$ for the resulting accuracy.

\subsubsection{The diffusion loss scheme}
The diffusion loss scheme proposed in \cite{r10,DiffusionLoss24} uses a DNN $u_\theta(t, x)$ that approximates the solution to the PDE \cref{pde} for $t \in [0, \tau]$, where $\tau$ is a stopping time.
The loss function is constructed using the matching of terminal and boundary conditions, as well as the integral form of the backward SDE, i.e. using the residual of the identity
\begin{equation}\label{eq_diffusion_loss_int}
    Y_{\tau} = Y_{0} + \int_0^\tau \phi(t, X_t, Y_t, Z_t, \cdots) dt + \int_0^\tau Z_t^T \sigma(t, X_t, Y_t, Z_t, \cdots) dW_t,
\end{equation}
where $Y_t = u_\theta(t, X_t)$ and $Z_t = \nabla u_\theta(t, X_t)$.
When $\tau = T = t_N$, with a na{\"i}ve numerical integration scheme under the same time discretization as in previous discussions, the component of the loss function brought by the above residual is equivalent to using $\lVert Y_N - \hat{Y}_N \rVert^2$ from the SDE matching approach of \cite{r9}.

It has also been discussed in \cite{r10} that as $T \to 0$, the limit of the diffusion loss differs with the loss of the PINN method by a constant multiplier.

\subsubsection{The deep martingale network}
The deep martingale network (DeepMartNet) \cite{r11,cai2024soc} uses a DNN $u_\theta(t, x)$ that approximates the solution to the PDE \cref{pde} for $t \in [0, T]$.
The idea is based on the martingale property of stochastic integral terms in \cref{eq_diffusion_loss_int}, i.e., for any stopping times $\tau_0 < \tau_1$,
\begin{equation}\label{eq_mart}
    Y_{\tau_1} - Y_{\tau_0} - \int_{\tau_0}^{\tau_1} \phi(t, X_t, Y_t, Z_t, \cdots) dt = \int_{\tau_0}^{\tau_1} Z_t^T \sigma(t, X_t, Y_t, Z_t, \cdots) dW_t 
\end{equation}
has zero expectation.
DeepMartNet samples segments of $[\tau_0, \tau_1]$ intervals that together cover $[0, T]$, on which massive trajectories are sampled to approximate the expectations.
When the functions $\mu$, $\phi$ and $\sigma$ do not depend on $z$ (and $Z_t$, accordingly), the method is gradient free, i.e. there's no need to explicitly evaluate $Z_t$ as the gradient of $u_\theta(t, X_t)$.
The trade-off is the efficiency on estimating the expectation, as the zero-mean stochastic integral indeed reduces the variance of the residual of \cref{eq_mart} under numerical discretizations.

\section{The shotgun method}\label{sec_method}
As discussed in the previous section, existing deep FBSDE methods simulate long trajectories of the backward SDE of $Y_t$, and utilize as much information collected along the trajectories as possible to form the loss functions for training DNNs.
In this section, we propose the ``shotgun method'' as an alternative approach that originates in the FBSDE formulation as well as the pioneer works of deep FBSDE methods, but also enjoys the flexibility of a PINN-type deep learning method.

\subsection{Single-step SDE discretization}

The derivation of the shotgun method starts from an observation on single-step simulations of discretized SDE trajectories using the Euler--Maruyama scheme.
Consider a time step size $0 < \Delta t \ll 1$.
Let $\eta = \Delta W_t \sim \mathcal{N}(0, \Delta t I_d)$ be a $d$-dimensional random vector that has the same distribution as the increment of a standard Brownian motion.
From the context of the PDE \cref{pde} and the associated FBSDE \cref{fbsde}, if we take time $t$ and $x = X_t$, and let $\mu$, $\sigma$ and $\phi$ be the values of the corresponding functions at $(t, X_t)$, respectively, then the single-step discretization of $X_{t+\Delta t}$ is given by
\begin{equation}\label{x_plus}
X_{t + \Delta t} \approx x + \mu \Delta t + \sigma \eta := x^+,
\end{equation}
and the single-step discretization of $Y_{t+\Delta t}$ is given by
\begin{equation*}
Y_{t + \Delta t} \approx Y_t + \phi \Delta t + Z_t^T \sigma \eta.
\end{equation*}
If, like in the FBSNNs method \cref{eq_fbsnn_scheme}, $Y_t$ and $Z_t$ are evaluated using a (smooth) test function $v(t, x)$ and its gradient $\nabla v(t, x)$, respectively, we accordingly define
\begin{equation}\label{y_plus}
y^+ = v(t, x) + \phi \Delta t + \nabla v(t, x)^T \sigma \eta.
\end{equation}

If $v(t, x) = u(t, x)$ is the solution to the PDE, then $y^+$ is a single-step approximation to $Y_{t+\Delta t} = u(t + \Delta t, X_{t+\Delta t})$.
In general, we claim that the approximation error partially consists of the residual of the PDE on the test function $v(t, x)$.
Namely, we write
\begin{equation}\label{v_plus}
\frac{v(t+\Delta t, x^+) - y^+}{\Delta t} = (\mathcal{L}[v](t, x) - \phi) + \varepsilon^+(t, x; v),
\end{equation}
where $\mathcal{L}$ is the operator defined in \cref{pde}, and $\varepsilon^+ = \varepsilon^+(t, x; v)$ is a random variable.
The size of $\varepsilon^+$ is considered to be small when $v$ is close to the solution to the PDE.
In fact, in the FBSNNs method \cref{eq_fbsnn_scheme}, if we take $(t, x) = (t_n, X_n)$, and $v = u_\theta$ as the DNN function, then the left-hand side of \cref{v_plus} is exactly $\Delta t^{-1} (Y_{n+1} - \hat{Y}_{n+1})$, whose square is a part of the loss function along the trajectory.

With the introduction of antithetic variates, i.e. with
\begin{equation}\label{x_minus}
x^- = x + \mu \Delta t - \sigma \eta, \quad
y^- = v(t, x) + \phi \Delta t - \nabla v(t, x)^T \sigma \eta,
\end{equation}
we also write
\begin{equation}\label{v_minus}
\frac{v(t+\Delta t, x^-) - y^-}{\Delta t} = (\mathcal{L}[v](t, x) - \phi) + \varepsilon^-(t, x; v),
\end{equation}
where $\varepsilon^- = \varepsilon^-(t, x; v)$ and $\varepsilon^+$ have the same distribution.
Thus, with the same instance of $\eta$, the average of \cref{v_plus} and \cref{v_minus} indicates
\begin{align}\label{eq_v_pm}
\begin{split}
   {}&{} \frac{1}{2} \frac{v(t+\Delta t, x^+) - y^+}{\Delta t} + \frac{1}{2} \frac{v(t+\Delta t, x^-) - y^-}{\Delta t} \\
   {}={} &\frac{v(t + \Delta t, x^+) + v(t + \Delta t, x^-) - 2v(t,x)}{2\Delta t} - \phi \\
   {} = {} & (\mathcal{L}[v](t, x) - \phi) + \varepsilon,
\end{split}
\end{align}
where
\begin{equation}\label{eq_eps}
\varepsilon=\varepsilon(t, x; v) = \frac{\varepsilon^+ + \varepsilon^-}{2}.
\end{equation}
We claim and prove in Appendix \ref{sect_appendix_proof} that as $\Delta t \to 0$,
\begin{equation}\label{eq_eps_mean_var}
\mathbb{E}[\varepsilon] = \mathcal{O}(\Delta t), \quad \mathrm{Var}(\varepsilon) = \mathcal{O}(1).
\end{equation}

The equation \cref{eq_v_pm} can be used as an estimator of the residual of the PDE on the test function $v$ at $(t, x)$, which is the term taken in the loss function of PINN.
It is worth pointing out that the gradient terms $\nabla v(t, x)$ from $y^+$ and $y^-$ cancel out in \cref{eq_v_pm}, so evaluating \cref{eq_v_pm} is almost always simpler when compared to \cref{v_plus}.
Since $\varepsilon^+$ and $\varepsilon^-$ have the same distribution, using the average of them is also guaranteed to have same or better error level in a probabilistic point of view.
Moreover, if the functions $\mu$, $\sigma$ and $\phi$ do not rely on the gradient of $v$, then \cref{eq_v_pm} is also gradient-free.

\subsection{The shotgun method}

The above discussion shows that a single-step time discretization \cref{v_plus} at $(t, x) = (t_n, X_n)$ in the deep BSDE methods is essentially an estimator on the residual of the PDE, which is used in the loss function of PINN.
In the perspective of PINN, the distribution of the collocation points to accumulate the loss function is required, but the collocation points forming complete stochastic process trajectories is often not necessary.
Thus, we consider a hybrid method using the distribution of the $X_t$ trajectories instead.
The method adopts the loss function of PINN, but the collocation points are taken from sparsely discretized trajectories which have approximately the same distribution as continuous trajectories, and the residual is estimated using \cref{eq_v_pm} instead of a straightforward evaluation of the Hessian matrix.

We discuss the following implementation details in the proposed method.
First, to maintain a uniform marginal distribution of time $t \in [0, T]$ on the sparsely discretized trajectories like in the limiting continuous version of the deep BSDE methods, we first set up a coarse uniform discretization of time into $N$ steps, with step size $\ell = T/N$.
Then, we sample a random variable $\delta \sim \mathrm{Uniform}(0, \ell)$, and take discretization points of $X_t$, $Y_t$ trajectories at
\begin{equation}
    t=0, \delta, \delta + \ell, \delta + 2\ell, \cdots, \delta + (N-1)\ell, T,
\end{equation}
respectively.

Second, as pointed out by \cref{eq_eps_mean_var}, the approximation error of the estimator \cref{eq_v_pm} has mean $\mathcal{O}(\Delta t)$ and variance $\mathcal{O}(1)$, so the variance may dominate the error.
Thus, we sample $M$ independent instances of $\eta \sim \mathcal{N}(0, \Delta t I_d)$ and take the average to reduce the variance of the approximation error.
The evaluation of the average is summarized in Algorithm \ref{algo_shotgun_step}.
Let the average of the approximation error $\varepsilon$ in \cref{eq_v_pm} from $M$ independent instances be denoted by $\epsilon^{(M)}$, then
\begin{equation}\label{eq_eps_M_mean_var}
    \mathbb{E}[\epsilon^{(M)}] = \mathcal{O}(\Delta t), \quad \mathrm{Var}(\epsilon^{(M)}) = \mathcal{O}(M^{-1}).
\end{equation}
In the rest of this paper, we call $\Delta t$ and $M$ the local time step size and the local batch size at each data point $(t, x)$, respectively, to distinguish with variables for trajectories and their discretization.
We emphasize that the local time step size $\Delta t$ may be independently chosen, and can be smaller than the coarse time step $\ell$ of the trajectories.

\begin{algorithm}[htbp]
\SetAlgoLined
\caption{Single-step residual estimate}\label{algo_shotgun_step}
\KwIn{PDE \cref{pde}, function $v(\cdot, \cdot)$, coordinates $(t, x)$, input parameters $y, z,\cdots$, local batch size $M$, local time step size $\Delta t$}
\KwOut{Estimated residual $r$}
$r \leftarrow 0$ \\
$\bar{x} \leftarrow x + \mu(t, x, y, z, \cdots) \Delta t$ \\
$\bar{\sigma} \leftarrow \sigma(t, x, y, z, \cdots)$ \\
\For{$m=1, \cdots, M$}
{
  Sample $\eta \sim \mathcal{N}(0, \sqrt{\Delta t} I_d)$ \\
  $x^+ \leftarrow \bar{x} + \bar{\sigma} \eta$ \\
  $x^- \leftarrow \bar{x} - \bar{\sigma} \eta$ \\
  $r \leftarrow r + \frac{1}{M}\left( \frac{v(t + \Delta t, x^+) + v(t + \Delta t, x^-) - 2 v(t, x)}{2 \Delta t} - \phi(t, x, y, z, \cdots) \right) $
}
\end{algorithm}

\begin{algorithm}[htbp]
\SetAlgoLined
\caption{Shotgun method}\label{algo_shotgun}
\KwIn{PDE \cref{pde}, distribution of initial spatial coordinate $p_0$, number of coarse trajectories $M_1$, number of time discretization steps $N$, local batch size $M$, local time step size $\Delta t$}
\KwOut{DNN function $u_\theta(t,x)$}
Initialize a DNN function $u_\theta(t, x)$ \\
$\ell \leftarrow T / N$ \\
\While{training}
{
$loss(\theta) \leftarrow 0$ \\
\For{$m=1,\cdots, M_1$}
{
Sample $\delta \sim \mathrm{Uniform}(0, \ell)$ \\
$t_0 \leftarrow 0$ \\
\For{$n=1,\cdots,N$}
{$t_n \leftarrow \delta + (n-1) \ell$}
$t_{N+1} \leftarrow T$\\
{Sample $X_0 \sim p_0$ \\
$Y_0 \leftarrow u_\theta(0, X_0)$ \\
\For{$n=0, \cdots, N$}
{
$Z_n \leftarrow \nabla u_\theta(t_n, X_n)$ \\
$\mu_n \leftarrow \mu(t_n, X_n, Y_n, Z_n, \cdots)$ \\
$\sigma_n \leftarrow \sigma(t_n, X_n, Y_n, Z_n, \cdots)$ \\
$\phi_n \leftarrow \phi(t_n, X_n, Y_n, Z_n, \cdots)$ \\
Sample $\eta_n \sim \mathcal{N}(0, (t_{n+1} - t_{n})I_d)$ \\
$X_{n+1} \leftarrow X_n + \mu_n (t_{n+1} - t_{n}) + \sigma_n \eta_n$ \\
$Y_{n+1} \leftarrow Y_n + \phi_n (t_{n+1} - t_{n}) + Z_n^T \sigma_n \eta_n$ \\
Apply Algorithm \ref{algo_shotgun_step} with $v \leftarrow u_\theta$ and $(t, x, y, z, \cdots) \leftarrow (t_n, X_n, Y_n, Z_n, \cdots)$ and collect residual $r$ \\
$loss(\theta) \leftarrow loss(\theta) + r^2 / (M_1 (N+1))$
}}
}
$loss(\theta) \leftarrow loss(\theta) + \lVert u_\theta(T, X_{N+1}) - g(X_{N+1}) \rVert^2 / M_1 $ \\
Apply stochastic gradient descent step to $loss(\theta)$
}
\end{algorithm}

The resulting method is a mixture of trajectory simulations like in deep FBSDE methods, and local PDE residual evaluation like in PINN.
We name the method a shotgun method, as it coarsely hunts for data points along the trajectories that in total covers the region of our interest, and takes multiple shots towards each single data point to reduce the error variance of the estimate.
An algorithm framework of the shotgun method is summarized in Algorithm \ref{algo_shotgun}.

Note that in the training of the shotgun method, the total number of collocation points is $(N+2)M$, and the backward propagation has $\mathcal{O}(N^2 M)$ computational cost.
Although it seemingly exceeds the $\mathcal{O}(N^2)$ cost of FBSNNs and the SDE matching method, it is later shown in the numerical tests in Section \ref{sec_result} that the shotgun method works with a much smaller $N$, as in previous deep BSDE methods a large $N = T / \Delta t$ is required for numerical precision, while the size of $M$ is not necessarily large.

\begin{remark}[Random finite difference scheme]
The residual estimation method in Algorithm \ref{algo_shotgun_step} can be interpreted as a finite difference scheme with randomness in steps.
Considering the practical accuracy of high-dimensional PDE solvers, using a finite difference scheme with moderately small $\Delta t$ is sufficient.
\cite{RAD20} and \cite{RFD23} use a randomized forward difference scheme for automatic differentiation in deep learning.
In \cite{r6}, the difference terms $(\partial_i v(t, x + \sigma \eta) - \partial_i v(t, x)) \eta_i$ are used to approximate the Hessian matrix of the function $v$ at $(t, x)$.
This approach requires a number of gradient calls at \emph{different} $(t, x+ \sigma \eta)$ positions.
\end{remark}

\subsection{Variants of the shotgun method}
Here we discuss a few variants of the shotgun method in Algorithm \ref{algo_shotgun}.

\begin{remark}[Gradient-free shotgun method]
If the gradients of $u(t, x)$ (such as $\nabla u(t, x)$) are not needed as a part of input variables of the functions in the PDE \cref{pde}, we may replace line 22 of Algorithm \ref{algo_shotgun} by
\begin{equation}
    Y_{n+1} \leftarrow u_\theta(t_{n+1}, X_{n+1})
\end{equation}
and omit line 16.
This leads to a gradient-free version of the algorithm.
\end{remark}

\begin{remark}[Time-independent PDEs]\label{remark_time_indep}
For PDEs that have a formulation close to \cref{pde} but do not have time terms, we may solve a time-dependent version together with an extra equation
\begin{equation}
    \partial_t u = 0.
\end{equation}
The extra equation is naturally satisfied when using a DNN function $u_\theta(x)$ that does not accept input of $t$.
\end{remark}

\begin{remark}[Shotgun-PINN method]\label{remark_shotgun_pinn}
We can use Algorithm \ref{algo_shotgun_step} as the \emph{backend} of the residual evaluation of differential equations in a PINN method.
In this paper we call such approach a shotgun-PINN method.
\end{remark}

\section{Numerical results}
\label{sec_result}

In this section, we illustrate numerical results of the proposed shotgun method in various test cases.
All numerical experiments, unless specifically stated, are performed on an NVIDIA A40 GPU with 48 GB of graphics memory.

\subsection{One-dimensional Laplace's equation}\label{sect_1d_laplace}
We begin with a 1-D example showing the influence of the local time step size $\Delta t$ and the local batch size $M$ to the resulting accuracy.

We study the 1-D Laplace's equation for $x \in [1, 5]$ with Dirichlet boundary conditions
\begin{align}\label{eq_Laplace1D}
    \begin{split}
        u''(x) &= f''(x), \quad x \in (1, 5), \\
        u(1) &= f(1), \\
        u(5) &= f(5).
    \end{split}
\end{align}
We take tests with $f(x) = f_j(x)$,  $j=1, 2, 3$, where
\begin{align}
    f_1(x)=x^2,\quad f_2(x)=\sin(x),\quad f_3(x)=e^{0.5x}.
\end{align}
The shotgun-PINN method (c.f. Remark \ref{remark_shotgun_pinn}) is applied to solve the above differential equation.
We use a fully-connected DNN $u_\theta(x)$ to approximate the solution $u(x)$, where $u_\theta$ is of size 1-64-64-1, i.e. $u_\theta$ has 3 layers, and the widths of the 2 hidden layers are 64, respectively.
The sine function $\sin x$ is used as the activation functions.
The loss function is given by
\begin{equation}
    L(\theta) = \sum_{i=1}^{100} r(x_i; \Delta t, M)^2 + |u_\theta(1) - f(1)|^2 + |u_\theta(5) - f(5)|^2,
\end{equation}
where $x_1, \cdots, x_{100} \sim \mathrm{Uniform}(1, 5)$ are i.i.d. sampled for each epoch, and $r(x_i; \Delta t, M)$ is the residual of the differential equation at $x_i$ evaluated using Algorithm \ref{algo_shotgun_step}.
The Adam optimizer is used to minimize $L(\theta)$ for 3000 epochs with learning rate $1 \times 10^{-3}$.
We uniformly sample 1000 independent points $y_1, \cdots, y_{1000} \sim \mathrm{Uniform}(1, 5)$ once training ends to estimate the $L^\infty$ error $\lVert u_\theta(x) - u(x) \rVert_{L^\infty(1, 5)}$ with
\begin{equation}\label{eq_numer1_error}
    \mathrm{error} = \max_{1 \le i \le 1000} |u_\theta(y_i) - u(y_i)|.
\end{equation}
For the case $f(x) = f_1(x)$, an error matrix with different values of local time step size $\Delta t$ and local batch size $M$ is illustrated in Table \ref{error_sin}.
The error is not much affected by $\Delta t$ when $\Delta t$ is sufficiently small as suggested by \cref{eq_eps_M_mean_var}.
This is indeed because $f_1(x)$ is quadratic, so that $\mathbb{E}[\varepsilon^{(M)}] = 0$ in this case.

\begin{figure}[H]
    \centering
    \begin{subfigure}[b]{0.32\textwidth}
        \includegraphics[width=\textwidth]{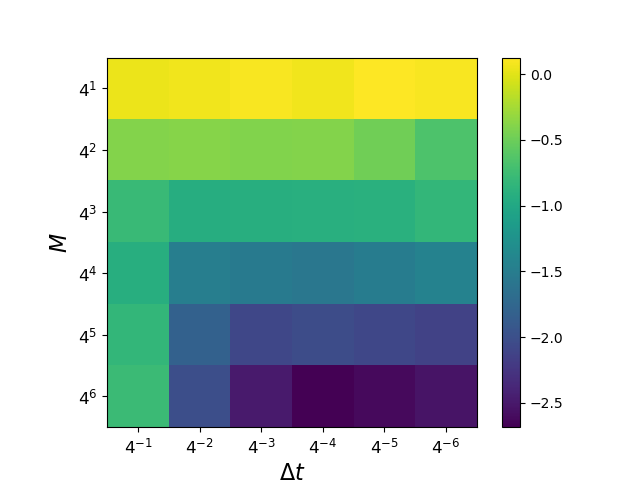}
        \caption{$f(x)=f_1(x) = x^2$}
        \label{fig:image1_x}
    \end{subfigure}
    \begin{subfigure}[b]{0.32\textwidth}
        \includegraphics[width=\textwidth]{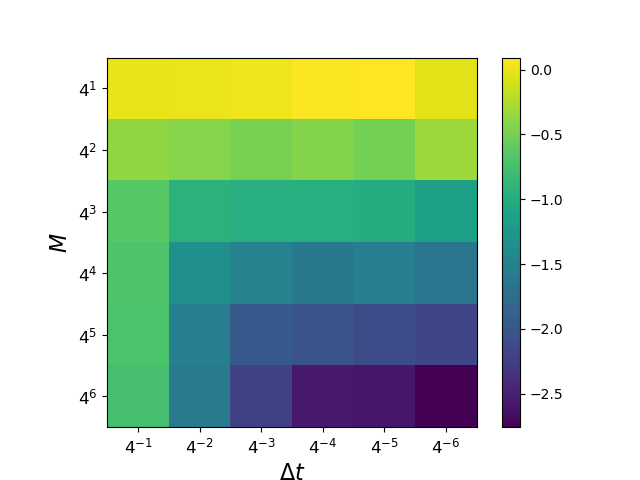}
        \caption{$f(x)=f_2(x)=\exp(0.5x)$}
        \label{fig:image1_exp}
    \end{subfigure}
    \begin{subfigure}[b]{0.32\textwidth}
        \includegraphics[width=\textwidth]{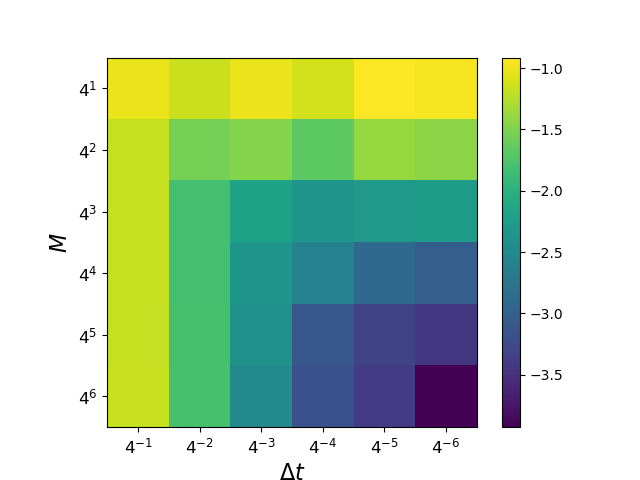}
        \caption{$f(x)=f_3(x)=\sin x$}
        \label{fig:image1_sin}
    \end{subfigure}
    \caption{Heat map of the error matrices of \cref{eq_numer1_error} (taken base-$10$ logarithm for illustration) for different test cases.}
    \label{fig:two_images_1D}
\end{figure}

For test cases $f(x) = f_2(x)$ and $f(x) = f_3(x)$, the errors are dominated by both the decrement of $\Delta t$ and the increment of $M$.
We compare the error matrices in heat maps in Figure \ref{fig:two_images_1D}.
It is visually clear that there are diagonal blocks that have a balance between errors dominated by $\Delta t$ and by $M$, respectively.

\begin{table}[H]
   \centering
   \begin{tabular}{|c|c|c|c|c|c|c|}
   \hline
   \diagbox{$M$}{error}{$\Delta t$}&$4^{-1}$&$4^{-2}$&$4^{-3}$&$4^{-4}$&$4^{-5}$&$4^{-6}$\\ 
    \hline
    $4^1$&1.12e-00&1.18e-00&1.26e-00&1.16e-00&1.34e-00&1.26e-00\\
    \hline
    $4^2$&4.00e-01&4.10e-01&3.91e-01&4.06e-01&3.28e-01&2.19e-01\\
    \hline
    $4^3$&1.65e-01&1.16e-01&1.19e-01&1.22e-01&1.25e-01&1.51e-01\\
    \hline
    $4^4$&1.19e-01&3.25e-02&2.88e-02&2.70e-02&3.15e-02&3.62e-02\\
    \hline
    $4^5$&1.49e-01&1.54e-02&8.15e-03&9.32e-03&8.12e-03&7.48e-03\\
    \hline
    $4^6$&1.72e-01&9.83e-03&3.30e-03&2.06e-03&2.45e-03&2.87e-03\\
    \hline
   \end{tabular}
   \caption{Error matrix for $f(x) = f_1(x)$ test case with different choices of local time step size $\Delta t$ and local batch size $M$.}\label{error_sin}
\end{table}

\subsection{100-dimensional Hamilton--Jacobi--Bellman equation}
\label{sect_hjb}
Originated in the linear-quadratic Gaussian (LQG) control problem, a special case of the Hamilton--Jacobi--Bellman (HJB) equation has been extensively tested numerically by various existing works \cite{r1,r2,cai2024soc,hu2023bias}, typically on a time interval $[0, T]$ with $T=1$.
In this section we use this test problem to demonstrate the capability of the shotgun method for solving the HJB equation on a longer time interval.

The 100-dimensional HJB equation is stated as follows.
Let $u(t, x)$ be defined on $[0, T]\times \mathbb{R}^d$ where $d = 100$, satisfying the parabolic PDE
\begin{align}\label{eq_hjb}
\begin{split}
    \partial_t u + \mathrm{Tr}[\nabla \nabla u] &= \lVert \nabla u \rVert^2, \\
    u(T, x) &= g(x),
\end{split}
\end{align}
where the squared norm $\lVert v \rVert^2 = v^T v$.
We let
\begin{equation}
    \mu(x) = 0, \quad \sigma(x) = \sqrt{2}I_d, \quad \phi(x, y, z) = \lVert z\rVert^2
\end{equation}
so that $X_t$ in the FBSDE formulation \cref{fbsde} is the controlled state process studied in this stochastic optimal control problem.
We set the initial state
\begin{equation}
X_0 = x_0 = (0, \cdots, 0).
\end{equation}
The terminal values at time $t=T$ is given with
\begin{equation}
    g(x) = \log (0.5(1 + \lVert x \rVert^2)).
\end{equation}
The exact solution to \cref{eq_hjb} is given by the formula
\begin{equation}\label{eq_hjb_exact}
    u(t, x) = -\log \mathbb{E} \left[ e^{-g(x + \sqrt{2}W_{T-t})} \right],
\end{equation}
whose evaluation requires Monte Carlo simulation.
In this paper, for each $(t, x)$ pair, we sample $10^5$ independent instances of $W_{T-t}$ and take the mean of the results as an estimate to the expectation in \cref{eq_hjb_exact}.
We believe the sample size is sufficient, as the sample standard deviation values are sufficiently small compared to the numerical errors of DNN solvers to be later illustrated.

We apply the shotgun method to solve the HJB equation \cref{eq_hjb} with time interval lengths $T=1$ and $T=10$, respectively.
The parameters used by the shotgun method are listed as follows.
The initial distribution of the stochastic process $X_t$ at $t=0$ is chosen as
\begin{equation}
    X_0 \sim \mathcal{N}(x_0, 2.5 \times 10^{-3} I).
\end{equation}
For the coarsely discretized trajectories, for each epoch, we pick $M_1 = 50$ trajectories, each discretized into $N=10$ steps when $T=1$, and into $N=100$ steps when $T=10$.
For the local parameters, we take $\Delta t = 4^{-5}$ and $M=64$.
We use a fully-connected DNN $u_\theta(t, x)$ to approximate the solution $u(t, x)$, where $u_\theta$ is of size ($d+1$)-512-512-512-512-1, i.e. the DNN has $5$ layers with $4$ hidden layers, and each hidden layer has width $512$.
The activation functions are the \texttt{mish} function given by
\begin{equation}
    \mathrm{mish}(x) = x \tanh \left( \log (1 + e^x) \right).
\end{equation}
The networks are trained using the Adam optimizer for $10000$ steps, with learning rate starting from $1\times 10^{-3}$ and decaying by a multiplier factor of $0.2$ per $2500$ steps.
Once the training ends, we simulate a total of $M_{test}=100$ SDE trajectories of $X_t$ with a fine uniform time discretization of $[0, T]$ with $N_{test}=100$ steps, and collect the relative error by the formula
\begin{equation}\label{eq_BSB_err}
    e_n(\omega) = \frac{|u_\theta(t_n, X_n(\omega)) - u(t_n, X_n(\omega))|}{|u(t_n, X_n(\omega))|}, \quad t_n = \frac{n}{N_{test}}T, \quad n=0,\cdots,N_{test},
\end{equation}
where $u(t, x)$ is the exact solution to the PDE, and each $\omega$ from the sample space refers to an instance of the Brownian motion $W_t$ in the FBSDE formulation.
We then collect the sample mean and the sample standard deviation of each $e_n$ using the results from these instances.

\begin{figure}[H]
    \centering
    \begin{subfigure}[b]{0.48\textwidth}
        \includegraphics[width=\textwidth]{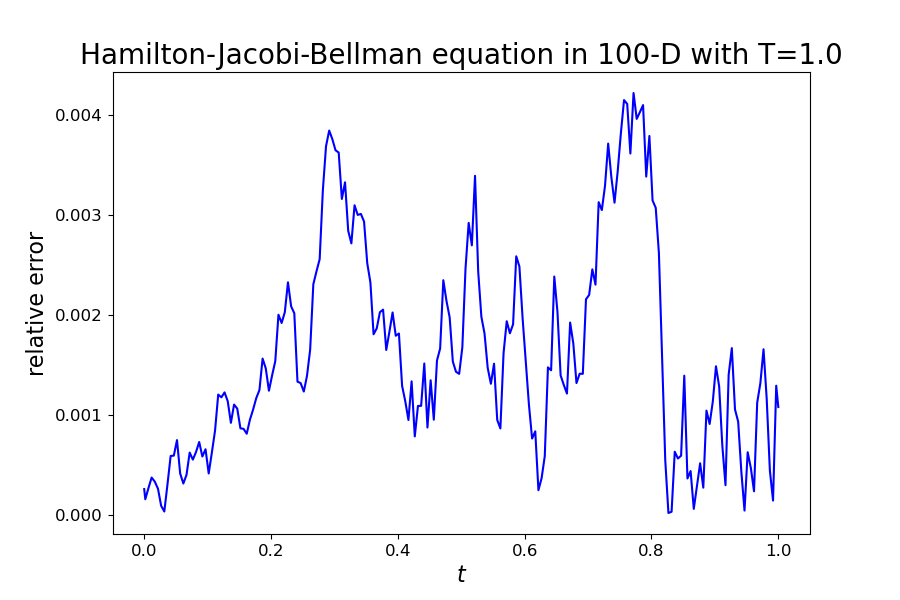}
        \caption{$T=1.0$}
        \label{fig:image2_HJB_100D_1}
    \end{subfigure}
    \begin{subfigure}[b]{0.48\textwidth}
        \includegraphics[width=\textwidth]{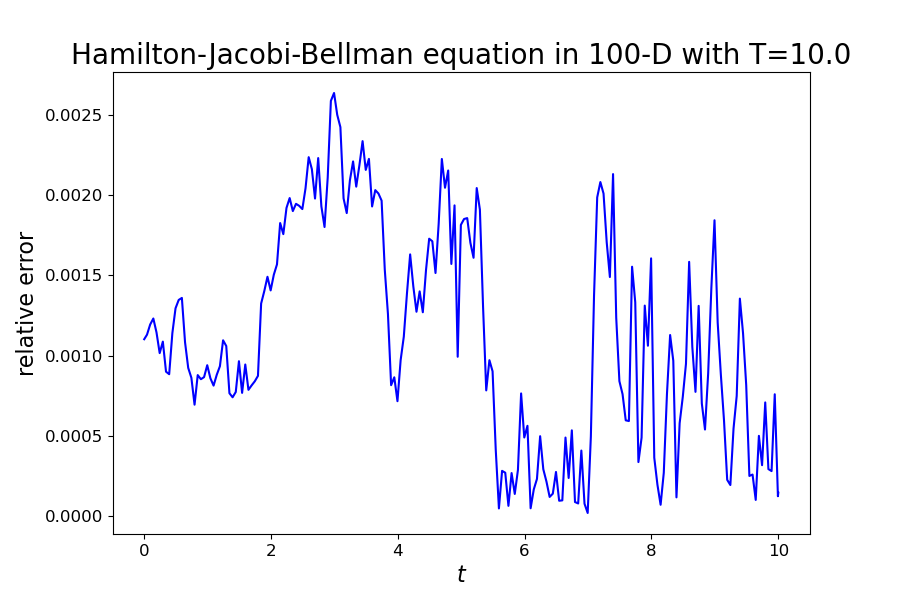}
        \caption{$T=10.0$}
        \label{fig:image2_HJB_100D_10}
    \end{subfigure}
    \caption{Relative errors $e_n(\omega)$ of the 100-dimensional Hamilton-Jacobi-Bellman equation \cref{eq_hjb} between model predictions and the reference solution along the time interval $[0, T]$ on a single realization of the underlying Brownian motion $W_t(\omega)$.}
    \label{fig:two_images_HJB_100D_error}
\end{figure}

\begin{figure}[H]
    \centering
    \begin{subfigure}[b]{0.48\textwidth}
        \includegraphics[width=\textwidth]{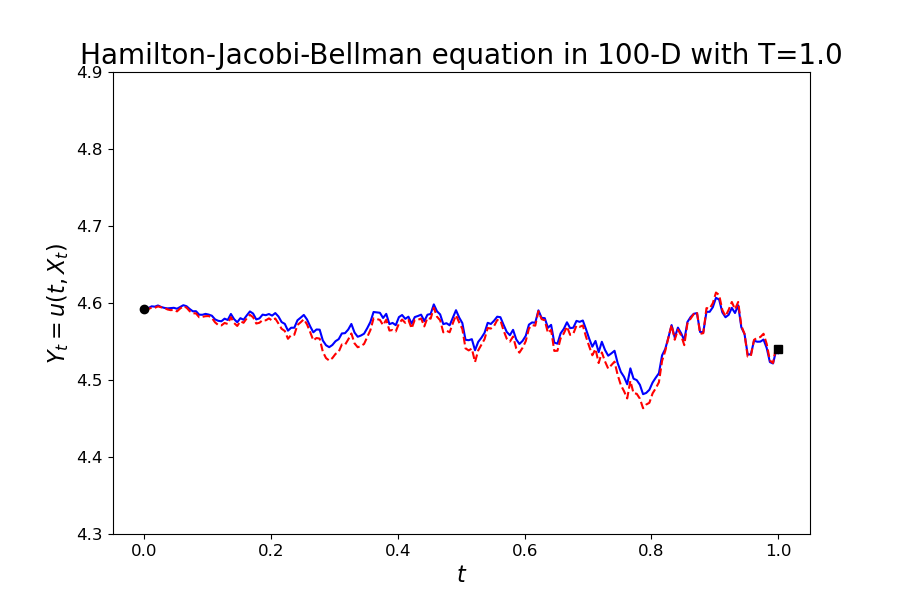}
        \caption{$T=1.0$}
        \label{fig:image1_HJB_100D_1}
    \end{subfigure}
    \begin{subfigure}[b]{0.48\textwidth}
        \includegraphics[width=\textwidth]{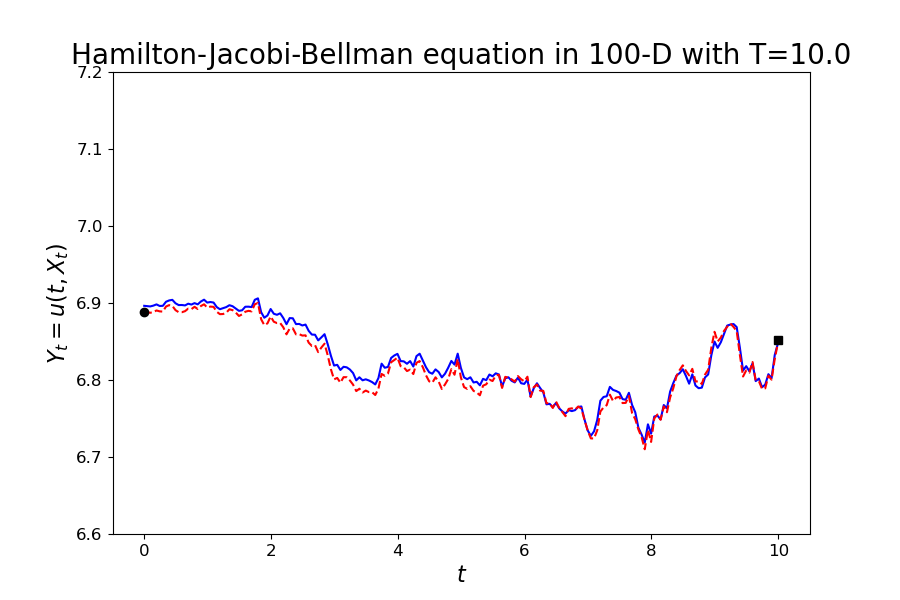}
        \caption{$T=10.0$}
        \label{fig:image1_HJB_100D_10}
    \end{subfigure}
    \caption{Predicting $Y_t$ trajectories in the FBSDE formulation using trained DNNs.
    The $X_t$ trajectories are simulated with fine time discretization on a realization of the underling Brownian motion $W_t(\omega)$ (same as in Figure \ref{fig:two_images_HJB_100D_error}).
    The red dashed lines are evaluated using a Monte Carlo simulation of $u(t, X_t)$ for reference, and the blue lines are predictions by the DNN $u_\theta(t, X_t)$.}
    \label{fig:two_images_HJB_100D_path}
\end{figure}

For both $T=1$ and $T=10$ cases, we plot the relative errors $e_n(\omega)$ and the predicted $Y_t(\omega)$ trajectories in Figure \ref{fig:two_images_HJB_100D_error} and Figure \ref{fig:two_images_HJB_100D_path}, respectively.
The shotgun method exhibits its capability of predicting $Y_t$ trajectories in these test cases.

\begin{remark}
The numerical challenge of the HJB equation in \cref{eq_hjb} can be (partially) seen in its FBSDE formulation.
By writing the SDE of $Y_t$ in the integral form
\begin{equation}
    Y_t = Y_0 + \int_0^t \lVert Z_s \rVert^2 ds + \int_0^t \sqrt{2} dW_t
\end{equation}
and taking expectation with $t=T$, one can look back at the PDE formulation and see that
\begin{equation}
    \mathbb{E}[g(\sqrt{2} W_T)] = u(0, x_0) + \mathbb{E} \left[\int_0^T \lVert \nabla u(s, \sqrt{2}W_s) \rVert^2 ds \right].
\end{equation}
In the above test cases, it is estimated by Monte Carlo methods that
\begin{equation}
\frac{\mathbb{E}[g(\sqrt{2} W_T)] - u(0, x_0)}{u(0, x_0)} < 0.01,
\end{equation}
so the gradient $\nabla u(s, \sqrt{2}W_s) = \nabla u(s, X_s)$ are relatively small along most of the trajectories, are hard to be captured by numerical solvers directly, and are even harder as the dimensionality increases.

\end{remark}

\subsection{100-dimensional Allen--Cahn equation}
We consider the Allen--Cahn equation from \cite{r1}, rewritten as a parabolic PDE of $u(t, x)$ with terminal condition:
\begin{align}\label{eq_allen_cahn}
\begin{split}
    \partial_t u + \Delta u + u - u^3 &= 0, \\
    u(T, x) &= g(x),
\end{split}
\end{align}
where $x \in \mathbb{R}^{d}$, $t \in [0, T]$ and $g(x) = 1/(2 + 0.4 \lVert x \rVert^2)$.
With $d = 100$ and $T = 0.3$, the initial value at
\begin{equation*}
x_0 = (0, \cdots, 0)
\end{equation*}
is approximately $u(0, x_0) \approx 0.0528$.
We let
\begin{equation}
    \mu(x) = 0, \quad \sigma(x) = \sqrt{2}I_d, \quad \phi(x, y) = -y + y^3
\end{equation}
in the associated FBSDE formulation, and solve the PDE \cref{eq_allen_cahn} using the shotgun method with the same setup as in the $T=0.3$ test case from the previous section.
The estimated initial value at $x_0$ from the trained DNN is
\begin{equation}
    u_\theta(0, x_0) = 0.05280 \cdots.
\end{equation}
\subsection{High-dimensional Black--Scholes--Barenblatt equation}
\label{sect_BSB}
In a previous work \cite{r9}, the SDE matching method having half-order accuracy in time discretization was proposed and numerically tested.
In this section, we compare the shotgun method proposed in this paper with the SDE matching method in \cite{r9} by numerically solving the Black--Scholes--Barenblatt (BSB) equation below of $u(t, x)$:
\begin{align}\label{eq_BSB}
\begin{split}
\partial_t u + \frac{1}{2} \mathrm{Tr}[\sigma^2 \mathrm{diag}(x x^T) \nabla \nabla u] &= (u - x^T \nabla u )r, \quad t \in [0, T], \quad x \in \mathbb{R}^{d} \\
u(T, x) &= \lVert x \rVert^2.
\end{split}
\end{align}
The exact solution of the above PDE is
\begin{equation}
    u(t, x) = e^{(r+\sigma^2)(T-t)}\lVert x \rVert^2.
\end{equation}
In our numerical tests, we set $\sigma = 0.4$ and $r=0.05$.
The initial value $u(0, x_0)$ at
\begin{equation}
x_0 = (1.0, 0.5, 1.0, 0.5, \cdots, 1.0, 0.5) \in \mathbb{R}^d
\end{equation}
is demanded in the BSB problem, i.e. with distribution density $\delta(x - x_0)$ of $X_0$ in the associated FBSDE.
In terms of the formulation used in the PDE \cref{pde}, in the FBSDE system \cref{fbsde} as well as in Algorithm \ref{algo_shotgun}, we use functions
\begin{equation}
    \mu(x) = 0, \quad \bar{\sigma}(x) = \sigma \mathrm{diag}(x), \quad \phi(x, y, z) = (y - x^T z)r.
\end{equation}

We first compare the shotgun method with the SDE matching method \cite{r9}, where the dimensionality $d = 100$ and the time interval length $T=1$.
In the SDE matching method, we choose $N=192$ from the example used in \cite{r9} as the number of time discretization steps with $\Delta t = T / N = 1/192$, and $M_1 = 50$ as the batch size.
In the shotgun method, for the coarsely discretized trajectories, for each epoch, we pick $M_1 = 50$ trajectories, each discretized into $N=10$ steps.
For the local parameters, we take $\Delta t = 4^{-5}$ and $M=4, 16$ and $64$, respectively.
Since the SDE of $X_t$ is given in the form $dX_t = \sigma \mathrm{diag}(X_t) dW_t$, the entries of $X_t$ are indeed geometric Brownian motions, and the value at $t=T$ greatly diverges as $T$ increases.
To address the numerical issue of matching the terminal condition, we apply the trick \cite{lu2021physics}
\begin{equation}\label{u_}
    u_\theta(t, x) = g(x) + (T-t) \Tilde{u}_\theta(t, x)
\end{equation}
to approximate the solution to the PDE \cref{eq_BSB} while ensure the matching of terminal values, where $\Tilde{u}_\theta$ is the actual DNN to be trained.
The DNN $\Tilde{u}_\theta$ is chosen as a fully-connected network of size ($d+1$)-512-512-512-512-1.
The activation functions are the \texttt{mish} function.
The networks are trained using the Adam optimizer for $10000$ steps, with learning rate starting from $1\times 10^{-3}$ and decaying by a multiplier factor of $0.2$ per $2500$ steps.
Once the training ends, the relative errors $e_n$ along SDE trajectories are collected similar to \cref{eq_BSB_err} in the previous test, as well as the estimates to its mean and its standard deviation.

\begin{figure}[htbp]
    \centering
    \begin{subfigure}[b]{0.48\textwidth}
        \includegraphics[width=\textwidth]{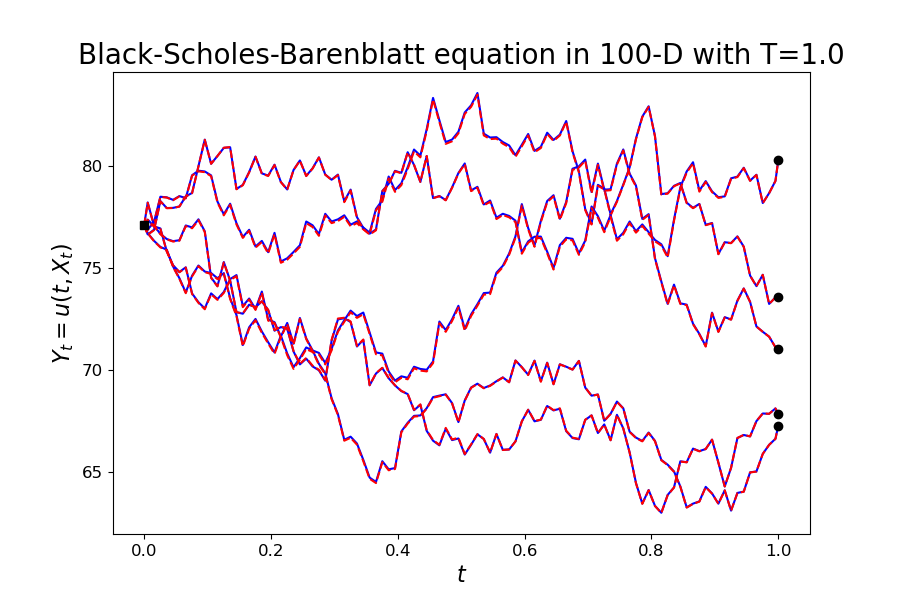}
        \caption{Shotgun, $M=64$}
        \label{fig:image1-1-SG100D-trajectory}
    \end{subfigure}
    \begin{subfigure}[b]{0.48\textwidth}
        \includegraphics[width=\textwidth]{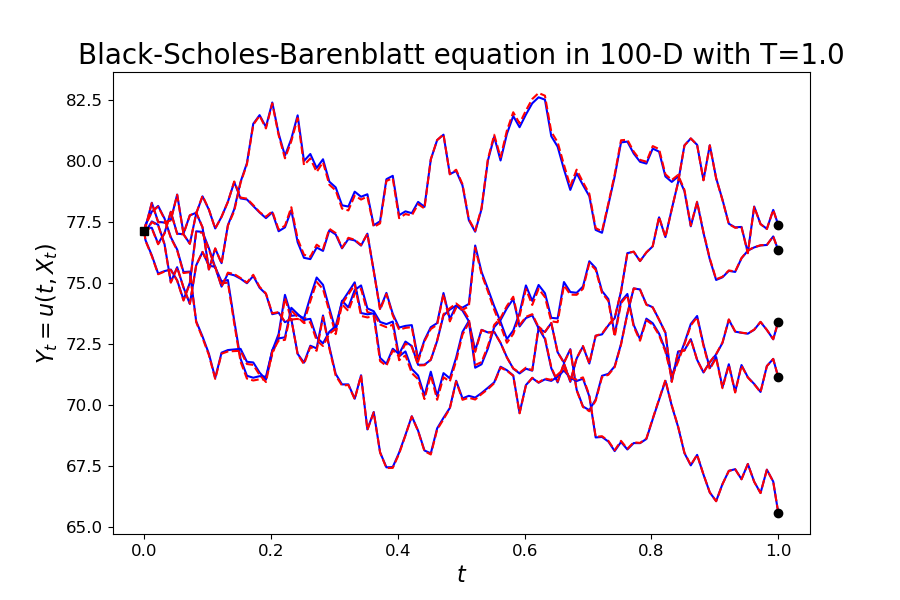}
        \caption{SDE matching, $N=192$}
        \label{fig:image1-1-FBSNN100D-trajectory}
    \end{subfigure}
    \caption{Predicting $Y_t$ trajectories in the FBSDE formulation of the 100-D Black-Scholes-Barenblatt equation \cref{eq_BSB} using trained DNNs.
    The $X_t$ trajectories are simulated with fine time discretization.
    The red dashed lines are evaluated using $u(t, X_t)$ for reference, and the blue lines are predictions of $u_\theta(t, X_t)$.}
    \label{fig:two_images-1-FBSNN and Shotgun100D trajectory}
\end{figure}

\begin{figure}[htbp]
    \centering
    \begin{subfigure}[b]{0.48\textwidth}
        \includegraphics[width=\textwidth]{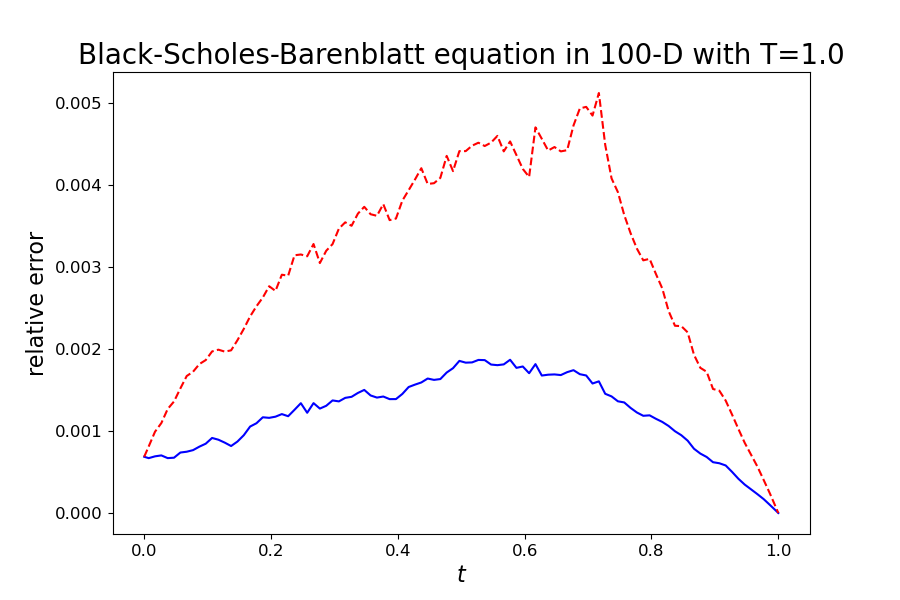}
       \caption{Shotgun, $M=4$}
        \label{fig:image2-1-SG100D-mean_2}
    \end{subfigure}
    \begin{subfigure}[b]{0.48\textwidth}
        \includegraphics[width=\textwidth]{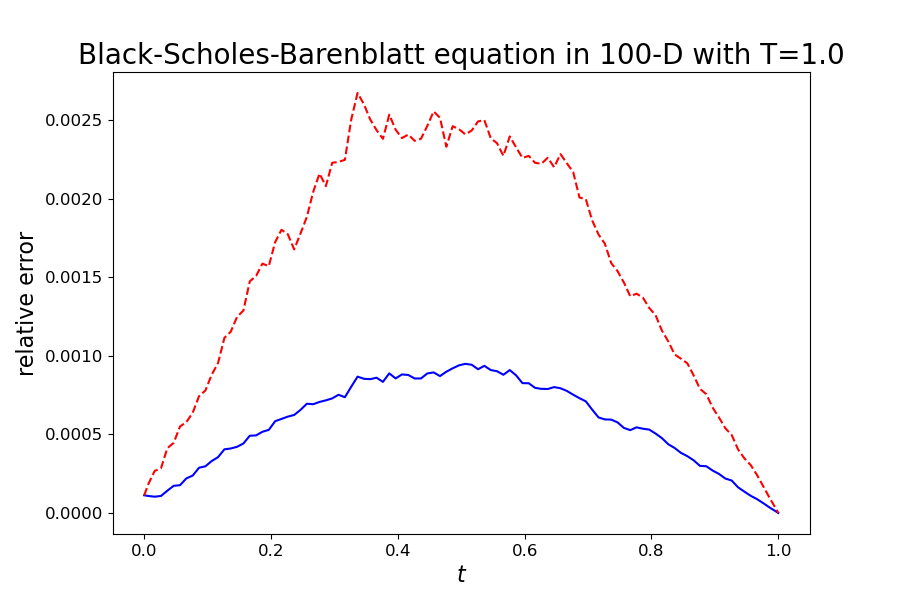}
       \caption{Shotgun, $M=16$}
        \label{fig:image2-1-SG100D-mean_4}
    \end{subfigure}
    \hfill
    \begin{subfigure}[b]{0.48\textwidth}
        \includegraphics[width=\textwidth]{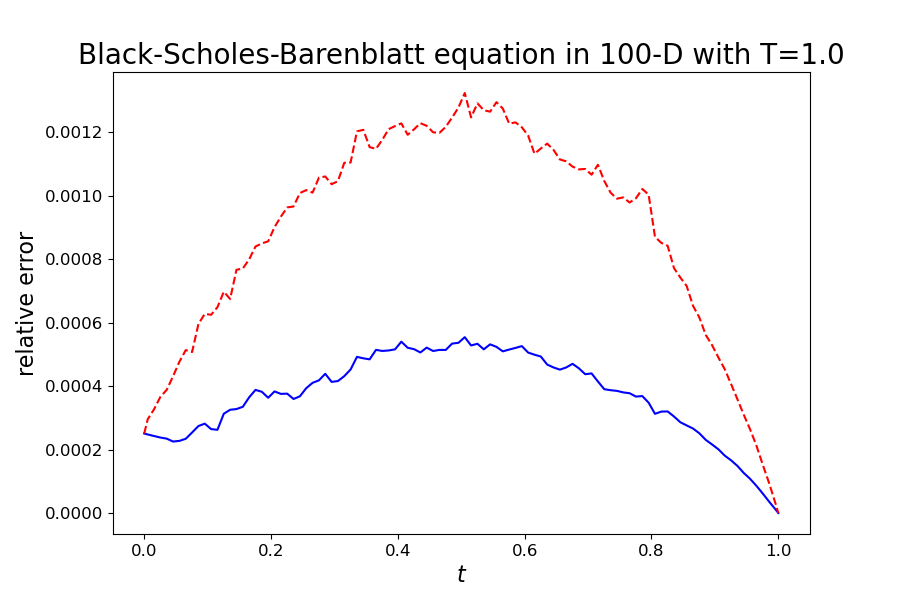}
       \caption{Shotgun, $M=64$}
        \label{fig:image2-1-SG100D-mean_6}
    \end{subfigure}
    \begin{subfigure}[b]{0.48\textwidth}
        \includegraphics[width=\textwidth]{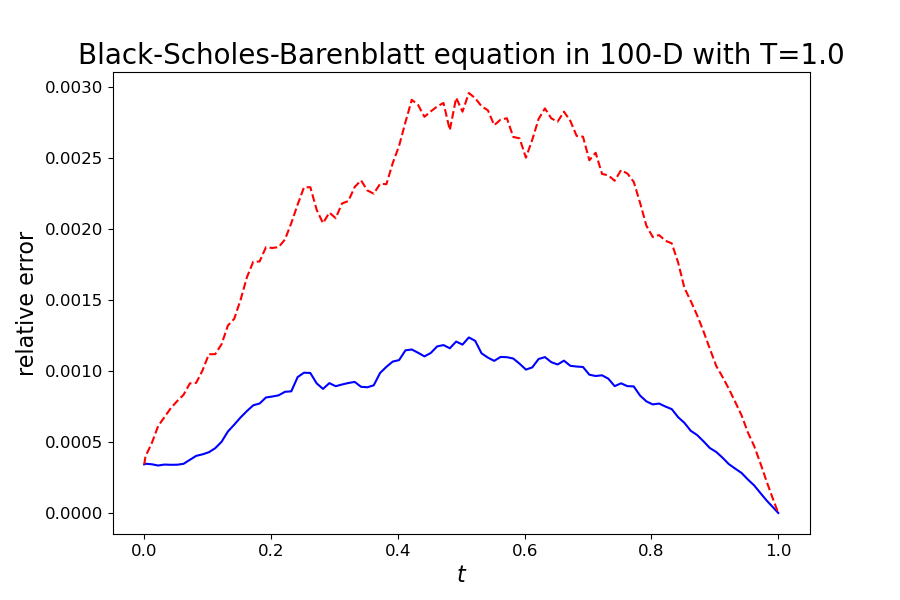}
       \caption{SDE matching, $N=192$}
        \label{fig:image2-1-FBSNN100D-mean_0}
    \end{subfigure}
    
    \caption{Test error $e_n$ of the 100-D Black-Scholes-Barenblatt equation \cref{eq_BSB} using the shotgun method with $M=4,16, 64$ and the SDE matching method \cite{r9}, respectively.
    The blue lines are estimates of $\mathbb{E}[e_n]$ along trajectories, and the red lines are estimates of $\mathbb{E}[e_n] + 2 \mathrm{std}(e_n)$.}
    \label{fig:two_images-1-FBSNN and Shotgun100D mean}
\end{figure}

Both the SDE matching method and the shotgun method are capable of predicting SDE trajectories, as shown in Figure \ref{fig:two_images-1-FBSNN and Shotgun100D trajectory}.
Meanwhile, the shotgun method exhibits superior precision when using local batch size $M=64$.
In Figure \ref{fig:two_images-1-FBSNN and Shotgun100D mean}, the test errors \cref{eq_BSB_err} are presented for using the shotgun method with $M=4$, $M=16$ and $M=64$, together with using the SDE matching method, respectively.
In the sequence of shotgun methods with increasing local batch size $M$, the error decays approximately at $\mathcal{O}(M^{-1/2})$ rate.
The shotgun method with $N=10$ and $M = 16$ has a close accuracy as the SDE matching method with $N=192$, where the total number of collocation points used for each trajectory is same, but the shotgun approach has lower cost in backward propagation.

Then, we further extend the test to dimensions $d > 100$ and with larger time interval size $T > 1$, respectively, which are rather challenging in performance for the SDE matching method.
In the test cases below, we apply local batch size $M=64$.
Because of the evident size variance of the $u(t, x)$ values, we solve the parabolic equation of $v(t, x)$
\begin{align}\label{eq_BSB_log}
\begin{split}
\partial_t v + \frac{1}{2} \mathrm{Tr}[\sigma^2 \mathrm{diag}(x x^T) \nabla \nabla v] &= r(1 - (\nabla v )^T x ) \\
&-\frac{1}{2}(\nabla v )^T \sigma^2 \mathrm{diag}(x x^T) \nabla v , \quad t \in [0, T], \quad x \in \mathbb{R}^{d}
\end{split}
\end{align}
instead of the original differential equation in \cref{eq_BSB}, where
\begin{equation}
    v(t, x) = \log u(t, x).
\end{equation}
The corresponding DNN to be trained is denoted by $\Tilde{v}_\theta(t, x)$, where
\begin{equation}
    v(t, x) \approx \log g(x) + (T - t) \Tilde{v}_\theta(t, x).
\end{equation}
Once the training of the DNN $\Tilde{v}_\theta$ ends, we examine the resulting accuracy on $u(t, x)$ of the original PDE \cref{eq_BSB} like in the previous comparison.

\subsubsection{Case $d > 100$}
\label{sect_case_big_d}
\begin{figure}[htbp]
    \centering
    \begin{subfigure}[b]{0.48\textwidth}
        \includegraphics[width=\textwidth]{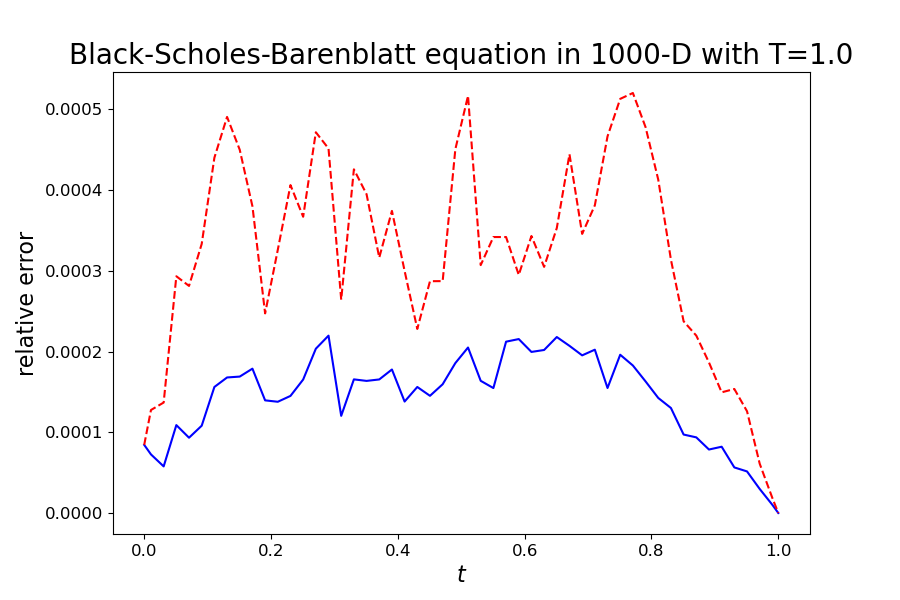}
       \caption{$d=1000$}
    \end{subfigure}
    \begin{subfigure}[b]{0.48\textwidth}
        \includegraphics[width=\textwidth]{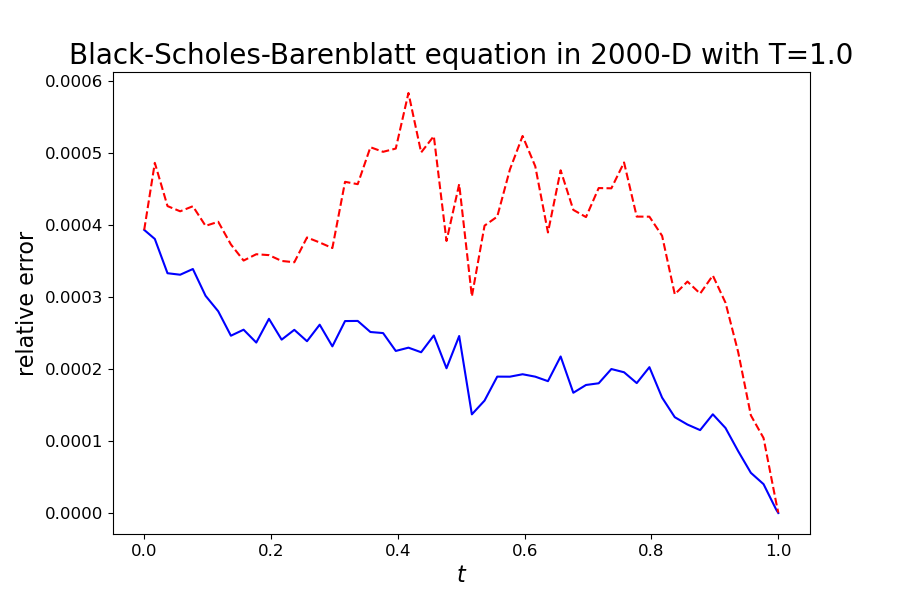}
       \caption{$d=2000$}
    \end{subfigure}
    \hfill
    \begin{subfigure}[b]{0.48\textwidth}
        \includegraphics[width=\textwidth]{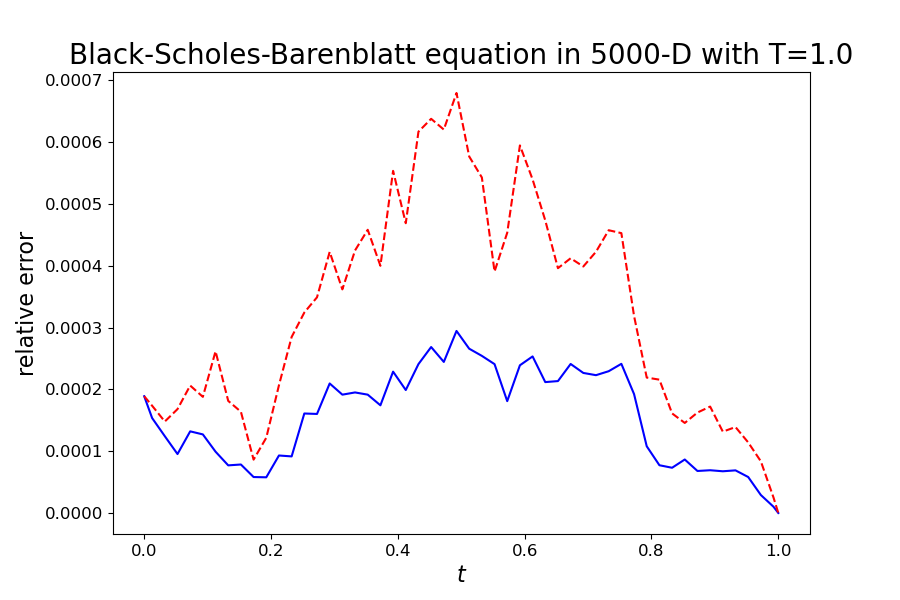}
       \caption{$d=5000$}
    \end{subfigure}
    \begin{subfigure}[b]{0.48\textwidth}
        \includegraphics[width=\textwidth]{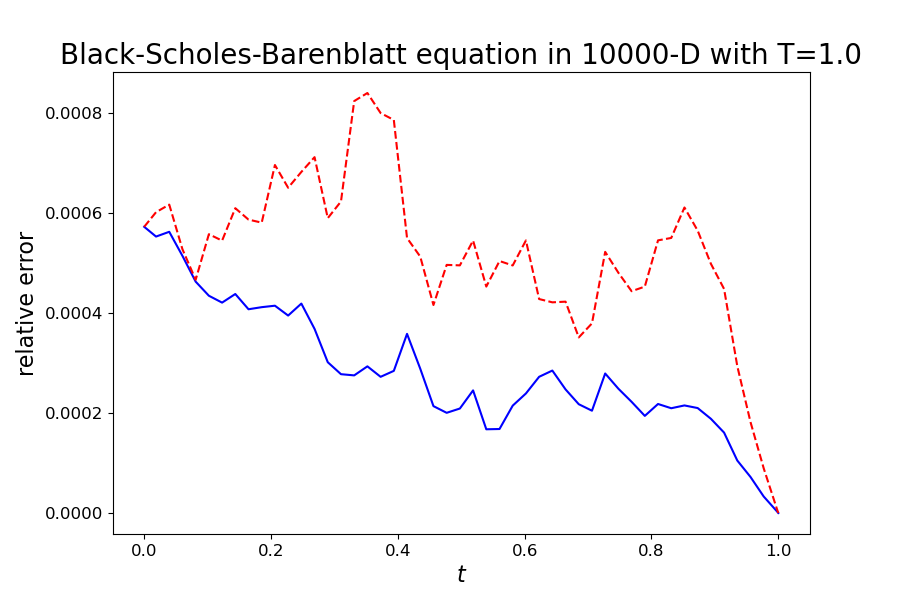}
       \caption{$d=10000$}
    \end{subfigure}
    
    \caption{Test error $e_n$ of the $d$-dimensional Black-Scholes-Barenblatt equation using the shotgun method for $d=1000, 2000, 5000$ and $10000$, respectively.
    The blue lines are estimates of $\mathbb{E}[e_n]$ along trajectories, and the red lines are estimates of $\mathbb{E}[e_n] + 2 \mathrm{std}(e_n)$.}
    \label{fig:logBSB_error_plots}
\end{figure}

Figure \ref{fig:logBSB_error_plots} shows the test errors $e_n$ in \cref{eq_BSB_err} with terminal time $T=1$ (and with trajectory discretization steps $N=10$) and with dimensionality $d=1000$, $2000$, $5000$ and $10000$\footnote{The test case of $d=10000$ was performed on an NVIDIA A100 GPU with 80 GB of graphics memory, so that test parameters including the batch sizes do not have to change.
A smaller size of hardware graphics memory is indeed sufficient if using a smaller batch size.}, respectively.

\begin{figure}[htbp]
    \centering
    \includegraphics[width=0.75\linewidth]{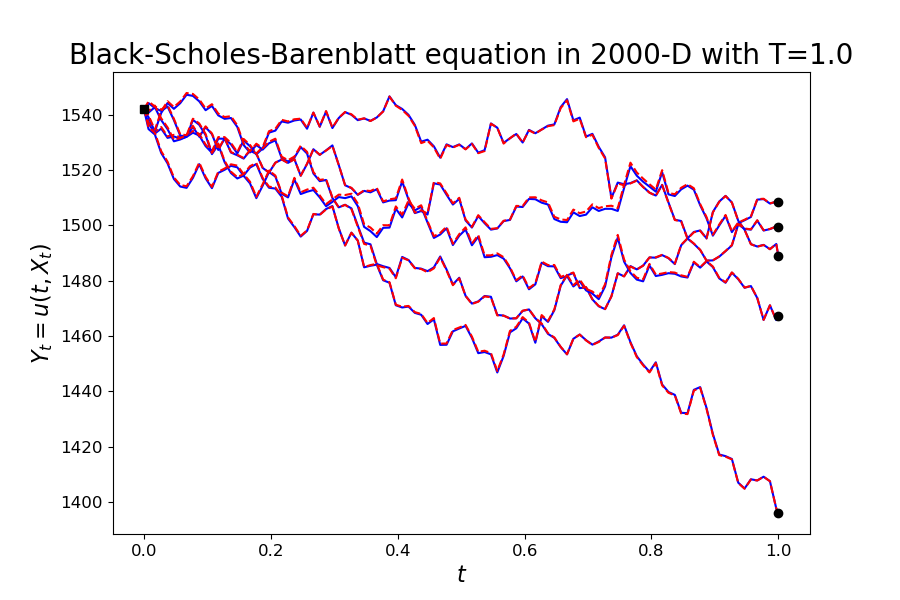}
    \caption{Predicting $Y_t$ trajectories in the FBSDE formulation of the $2000$-dimensional BSB equation \cref{eq_BSB} using trained DNNs.
    The red dashed lines are evaluated using $u(t, X_t) = \exp(v(t, X_t))$ for reference, and the blue lines are predictions using trained $\Tilde{v}_\theta (t, X_t )$.}
    \label{fig:BSB_2000}
\end{figure}

We also plot a few trajectories of $Y_t$ in the FBSDE formulation for illustration in Figure \ref{fig:BSB_2000}, where $d=2000$.
The shotgun method achieves competitive resulting accuracy in terms of relative error for these high-dimensional test cases.

\subsubsection{Case $T=10$}

\begin{figure}[htbp]
    \centering
    \begin{subfigure}[b]{0.48\textwidth}
        \includegraphics[width=\textwidth]{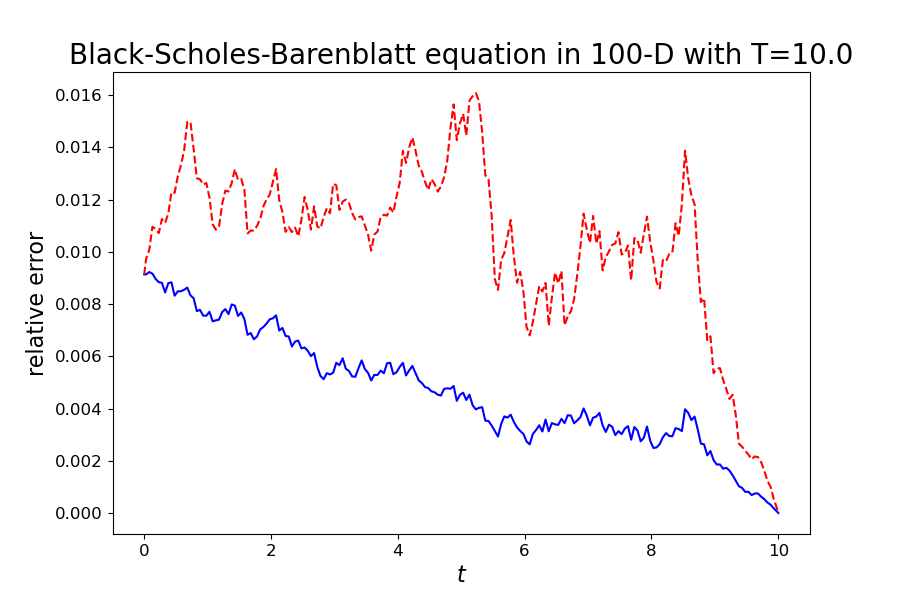}
       \caption{$d=100$}
    \end{subfigure}
    \begin{subfigure}[b]{0.48\textwidth}
        \includegraphics[width=\textwidth]{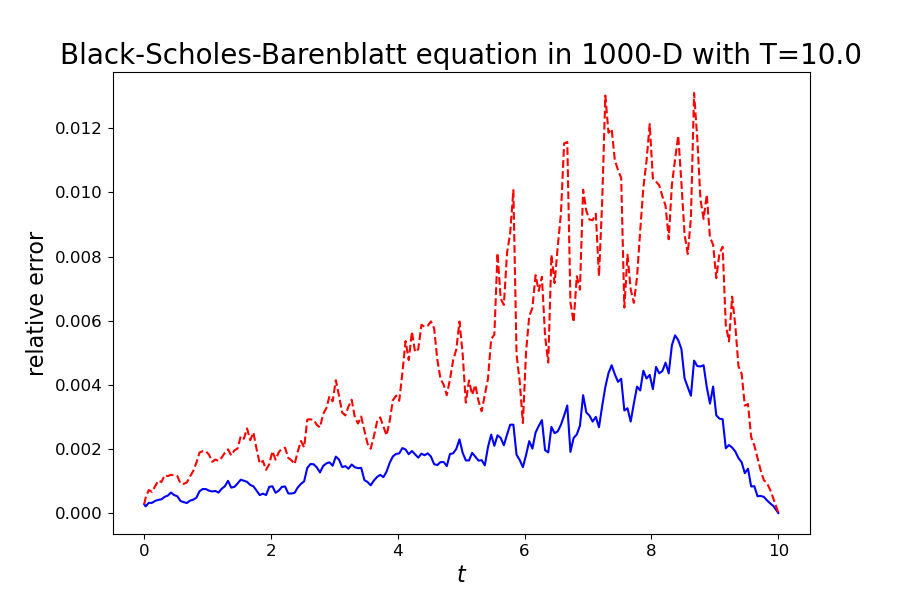}
       \caption{$d=1000$}
    \end{subfigure}
    \caption{Test error $e_n$ of the $d$-dimensional Black-Scholes-Barenblatt equation with $T=10$ using the shotgun method for $d=100$ and $d=1000$, respectively.
    The blue lines are estimates of $\mathbb{E}[e_n]$ along trajectories, and the red lines are estimates of $\mathbb{E}[e_n] + 2 \mathrm{std}(e_n)$.}
    \label{fig:BSB_T10_error}
\end{figure}

\begin{figure}[htbp]
    \centering
    \begin{subfigure}[b]{0.48\textwidth}
        \includegraphics[width=\textwidth]{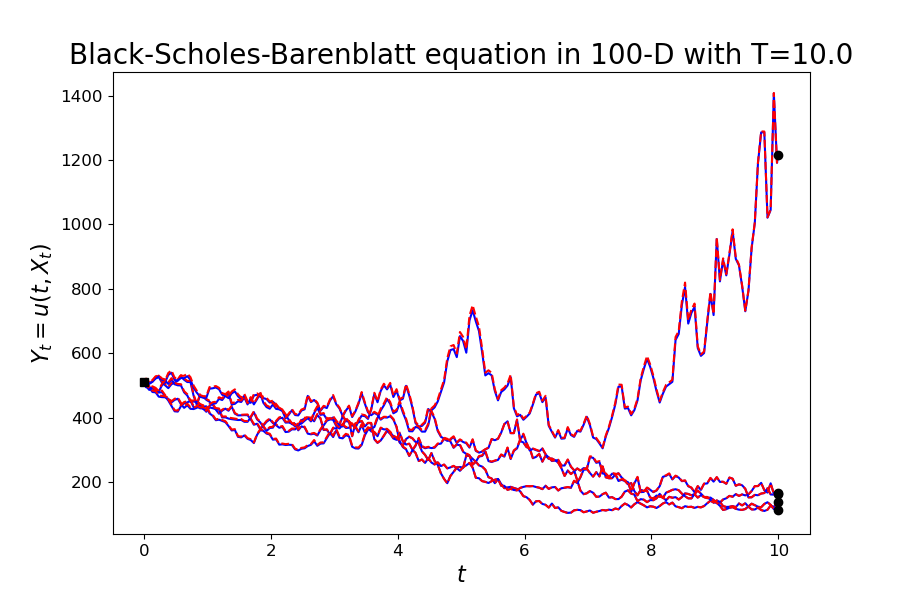}
        \caption{$d=100$}
    \end{subfigure}
    \begin{subfigure}[b]{0.48\textwidth}
        \includegraphics[width=\textwidth]{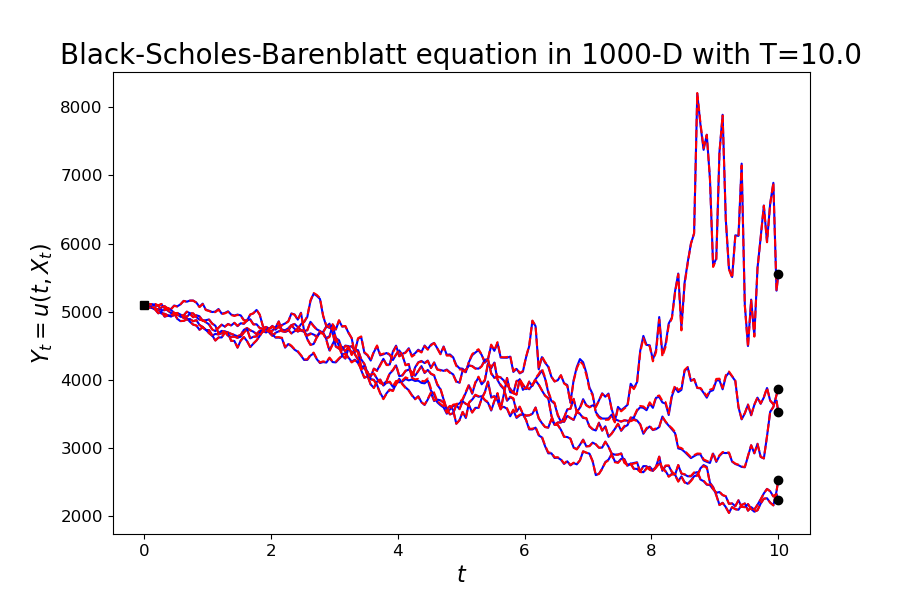}
        \caption{$d=1000$}
    \end{subfigure}
    \caption{Predicting $Y_t$ trajectories in the FBSDE formulation of the Black-Scholes-Barenblatt equation \cref{eq_BSB} with $T=10$ using trained DNNs.
    The red dashed lines are evaluated using $u(t, X_t)$ for reference, and the blue lines are predictions of $u_\theta(t, X_t)$.}
    \label{fig:BSB_T10_path}
\end{figure}

Figure \ref{fig:BSB_T10_error} shows the test errors $e_n$ (see \cref{eq_BSB_err}) with terminal time $T=10$ (and with trajectory discretization steps $N=100$) and with dimensionality $d=100$ and $1000$, respectively.
With $T=10$, the values of $X_t$ and $Y_t$ vary a lot more along trajectories compared to the $T=1$ cases, as presented in the $Y_t$ trajectory plots in Figure \ref{fig:BSB_T10_path}.
This causes a reduction in the resulting precision, but the overall relative error is still within 2\%.

\section{Discussion}
\label{sec_discussion}

\subsection{Benefits of the shotgun method}
As stated in Section \ref{sec_method}, the shotgun method is a hybrid method that takes advantages of deep FBSDE methods and the PINN method.
In the itemized comparison below, the shotgun method always picks the merit between deep FBSDE methods and PINN in each aspect.
\begin{enumerate}
    \item Data distribution.
    In high-dimensional parabolic PDEs of our concern, the SDE trajectories are preferred, meanwhile a vanilla PINN method has no clue about where to sample.
    The shotgun method effectively utilizes the distribution of SDE trajectories.
    \item Accuracy.
    When using a small local step size $\Delta t$ in the shotgun method, the resulting accuracy reaches a level that vanilla PINN usually achieves in practical tests, and outperforms those from time-consuming deep FBSDE experiments with small time step size such as in \cite{r9}.
    Note that tuning $\Delta t$ is independent from the performance of shotgun method, but has immediate impact on those of deep FBSDE methods, due to the evaluation on full trajectories.
    \item Training complexity and capability of solving long-time problems.
    In some existing deep FBSDE methods, the DNN output values are recursively used along SDE trajectories, causing increasing complexity in backward propagation as the recursion goes on.
    In PINN, the residuals at data points are independently evaluated.
    However, the Hessian matrices are very expensive for high-dimensional problems.
    The shotgun method avoids both problems.
    For long-time problems, coarse trajectories are used instead of fine discretization for efficiency.
\end{enumerate}

\subsection{Why estimators with big variance work?}
\label{sect_big_variance}
Intuitively, the error $\varepsilon$ of the residual estimator in \cref{eq_v_pm} has $\mathcal{O}(1)$ variance, so it brings $\mathcal{O}(1)$ error to the residuals.
However, the numerical tests in Section \ref{sec_result} indicate a different behavior.
To explain the convergence with the existence of apparent error variance, especially in the test case studied in Section \ref{sect_case_big_d} with
\begin{equation*}
    M = 64 \ll d = 10000,
\end{equation*}
we study the following simplified problem with a heuristic (yet non-rigorous) derivation.

Let $\theta$ be the vector of trainable variables of a DNN to be trained.
Let
\begin{equation}
L(\theta) = \sum_{k=1}^{K} f_k(\theta)^2
\end{equation}
be the reference loss function to be minimized with steepest descent method, where $f_1(\theta), \cdots, f_{K}(\theta)$ are scalar functions.
On the $n$-th step, the $f_k(\theta)$ residuals are perturbed with random functions that are independent from past iterations.
Namely, on the $n$-th step the modified loss function
\begin{equation}
L_n(\theta) = \sum_{k=1}^{K} \left(f_k(\theta) + \varepsilon_k^{(n)}(\theta)\right)^2
\end{equation}
is minimized instead.
Let $\theta_n$ be the value of $\theta$ in the iteration before step $n$.
Let $\alpha > 0$ be the learning rate.
Then,
\begin{equation}
\theta_{n+1} = \theta_n - \alpha \nabla L_n(\theta_n).
\end{equation}
When the training is close to convergence, we make assumptions that each $f_k(\theta)$ is a nontrivial \emph{linear} function of $\theta$, and that each $\varepsilon_k^{(n)}(\theta) = \varepsilon_k^{(n)}$ does \emph{not} depend on $\theta$ (i.e. the influence of existing perturbation to the network is negligible).
Then, for each $f_k$ function,
\begin{align}
\begin{split}
f_k(\theta_{n+1}) &= f_k(\theta_n) + \nabla f_k^T (\theta_{n+1} - \theta_n) \\
&= f_k(\theta_n) - 2 \alpha \nabla f_k^T \sum_{l=1}^{K} f_l(\theta_n) \nabla f_l - 2 \alpha \nabla f_k^T \sum_{l=1}^{K} \varepsilon_l^{(n)}\nabla f_l.
\end{split}
\end{align}
In a matrix-vector form, let
\begin{equation}
f(\theta) = \begin{bmatrix}
f_1(\theta) & \cdots & f_K(\theta)
\end{bmatrix}^T, \quad \varepsilon^{(n)} = \begin{bmatrix}
\varepsilon^{(n)}_1 & \cdots & \varepsilon^{(n)}_K
\end{bmatrix}^T, \quad D = [\nabla f_k \cdot \nabla f_l]_{k,l},
\end{equation}
we then have
\begin{equation}
f(\theta_{n+1}) = \left( I - 2 \alpha D \right) f(\theta_{n}) - 2 \alpha D \varepsilon^{(n)}.
\end{equation}
Suppose
\begin{equation}
\mathbb{E}[\varepsilon^{(n)}] = 0, \quad \mathrm{Cov}(\varepsilon^{(n)}) = V
\end{equation}
for simplicity.
Let $F_n = \mathrm{Cov}(f(\theta_n))$, then there holds the recursion
\begin{equation}
F_{n+1} = (I - 2 \alpha D) F_{n} (I - 2 \alpha D) + 4 \alpha^2 DVD.
\end{equation}
When $\alpha$ is sufficiently small, the limit of $F_n$ as $n \to \infty$ has order $\mathcal{O}(\alpha V)$, i.e., a sufficiently small learning rate $\alpha$ compensates for a potentially large variance $V$ in the resulting variance of the residuals.

\subsection{Choice of local batch size}

It follows from the numerical tests and the previous discussion that the local batch size $M$ can be much smaller than the dimensionality $d$.
We recommend adopting
\begin{equation}
    1 < M \ll d
\end{equation}
as a general rule.
Note that a deterministic central difference scheme for the second-order derivatives in the PDE \cref{pde} takes $\mathcal{O}(d)$ terms, so it is not worthy letting $M$ be as large as the size of $d$.

\subsection{Scopes of PDEs and limitations}

For the Black--Scholes--Barenblatt equation, we emphasize that the log-BSB version \cref{eq_BSB_log} is \emph{nonlinear} meanwhile the original version \cref{eq_BSB} is linear.
Numerical results in Section \ref{sect_BSB} indicate that the shotgun method is capable of solving nonlinear PDEs.

For the Hamilton--Jacobi--Bellman equation, it has been pointed out in \cite{wang20222} that the squared $L^2$ residual is not an ideal choice of PINN, and the authors recommend using $L^\infty$ norm instead.
Alternatively, \cite{cai2024soc} and \cite{sdgd2024} apply adversarial training.
Despite having seemingly good accuracy in Section \ref{sect_hjb}, the shotgun method can not get rid of the above limitation, in the sense that the method requires squared $L^2$ loss form as discussed in Section \ref{sect_big_variance}, and it will not benefit from an $L^\infty$ loss form.
\section{Conclusion}
\label{sec_conclusion}
We proposed a deep shotgun method that solves parabolic PDEs in high-dimensional spaces.
The shotgun method is a hybrid method that lies between deep FBSDE methods and the PINN method, enjoying the benefit from both sides.
The numerical tests boost the scale of PDE problems by $100$ times compared to our previous work in terms of the pair $(d, T)$ of dimensionality and time interval length.
We leave the study of random finite difference schemes to a future work.

\begin{acknowledgements}
The author of Wenzhong Zhang acknowledges financial support by National Science Foundation of China (Grant No. 92270205, No. 1220012530) and by National Key R\&D Program of China (No. 2022YFA1005202, No. 2022YFA1005203).
\end{acknowledgements}

%
\section*{Declarations}

Data will be made available upon request.

\section*{CRediT authorship contribution statements}

Wenjun Xu: Conceptualization, Methodology, Software, Validation and Writing – original draft . Wenzhong Zhang: Conceptualization, Investigation, Methodology, Supervision, Writing – original draft, Writing – review and editing.
\section*{Conflict of interest}
The authors declare that they have no conflict of interest.

\appendix
\section{Derivation of eq. (\ref{eq_eps_mean_var})}\label{sect_appendix_proof}

Here we assume the function $v(t, x)$ is sufficiently smooth with converging Taylor series.
We also assume that $ \lVert \eta \rVert \ll 1$, which has high probability when $\Delta t \ll 1$.

In the third-order Taylor series of $\epsilon \Delta t$ at $(t, x)$, terms with an odd order of $\eta$ are canceled out.
It follows that
\begin{align}
\begin{split}
    \varepsilon ={}&{} \left( \frac{1}{2}\partial_{tt}v + \frac{1}{2} \mu^T \nabla\nabla v \mu + \mu^T \nabla \partial_t v \right) \Delta t^2 + \\
    {}&{}+ \left( \frac{1}{2} \eta^T \sigma^T \nabla \nabla v \sigma \eta - \frac{1}{2} \mathrm{Tr}[\sigma \sigma^T \nabla \nabla v] \Delta t \right) \\
    {}&{} + \left( \frac{1}{6} \partial_{ttt} v + \frac{1}{2}\mu^T \nabla \partial_{tt} v + \frac{1}{2} \mu^T \nabla \nabla \partial_t v \mu \right) \Delta t^3 \\
    {} & {} + \frac{1}{2} \eta^T \sigma^T \nabla \nabla \partial_t v \sigma \eta \Delta t \\
    {}&{} + tail_4 \\
    :={}&{} e_1 + e_2 + e_3 + e_4 + tail_4.
\end{split}
\end{align}
The terms $e_1$ and $e_3$ are deterministic.
For the term $e_2$, let $\delta = \Delta t^{-1/2} \eta$, and $H = \sigma^T \nabla \nabla v \sigma$.
Diagonalize the real matrix $H$ with an orthogonal matrix $U$ by
\begin{equation}
H = U^T \Sigma_H U, \quad \Sigma_H = \mathrm{diag}(\sigma_H^1, \cdots, \sigma_H^d), \quad U^T = U^{-1}.
\end{equation}
Then, $U \delta \sim \mathcal{N}(0, I_d)$, and one can easily tell that
\begin{equation}
\mathbb{E}\left[ \frac{\Delta t}{2} \delta^T H \delta \right] = \frac{\Delta t}{2} \mathrm{Tr}[H], \quad
\mathrm{Var}\left( \frac{\Delta t}{2} \delta^T H \delta \right) = \frac{\Delta t^2}{2} \mathrm{Tr}[H^2],
\end{equation}
so
\begin{equation}
    \mathbb{E}[e_2] = 0, \quad \mathrm{Var}(e_2) = \frac{1}{2} \mathrm{Tr}[(\sigma^T \nabla \nabla v \sigma)^2] \Delta t^2.
\end{equation}
Similarly for $e_4$,
\begin{equation}
    \mathbb{E}[e_4] = \frac{1}{2} \mathrm{Tr} [\sigma^T \nabla \nabla \partial_t v \sigma] \Delta t^2,\quad 
    \mathrm{Var}(e_4) = \frac{1}{2} \mathrm{Tr} [(\sigma^T \nabla \nabla \partial_t v \sigma)^2] \Delta t^4.
\end{equation}
For the $tail_4$ term, it is straightforward that
\begin{equation}
    tail_4 = \mathcal{O}(\Delta t^4 + \lVert x^+ - x \rVert^4 + \lVert x^- - x \rVert^4 ) = \mathcal{O}( \Delta t^4 + \lVert \mu \Delta t \rVert^4 + \lVert \sigma \eta \rVert^4),
\end{equation}
so
\begin{equation}
    \mathbb{E}[tail_4] = \mathcal{O}(\Delta t^2), \quad \mathrm{Var}(tail_4) = \mathcal{O}(\Delta t^4).
\end{equation}
In summary,
\begin{equation}
    \mathbb{E}[\varepsilon \Delta t] = \mathbb{E}[e_1] + \cdots + \mathbb{E}[e_4] + \mathbb{E}[tail_4] = \mathcal{O}(\Delta t^2),
\end{equation}
and
\begin{equation}
    \mathrm{Var}(\varepsilon \Delta t) \le 5 \left( \mathrm{Var}(e_1) + \cdots + \mathrm{Var}(e_4) + \mathrm{Var}(tail_4) \right) = \mathcal{O}(\Delta t^2).
\end{equation}
As $\Delta t \to 0$, the variance of $\varepsilon \Delta t$ is dominated by the contribution of the term $e_2$.

\bibliographystyle{spmpsci}      
\bibliography{references.bib}

\end{document}